\documentclass[reqno,12pt]{amsart}\usepackage{latexsym,amssymb,amsthm,amsmath,amscd,a4wide}

\theoremstyle{plain}
\newtheorem{theorem}{Theorem}[section]
\newtheorem{lemma}{Lemma}[section]
\newtheorem{proposition}{Proposition}[section]
\newtheorem{corollary}{Corollary}[section]

\theoremstyle{remark}
\newtheorem{remark}{Remark}[section]
\theoremstyle{remark}

\theoremstyle{definition}
\newtheorem{definition}{Definition}[section]
\newtheorem{example}{Example}[section]
\numberwithin{equation}{section}

\def\<{\left< }
\def\>{\right> }
\def\E4{\mathbb E^4 }

\def\E4{\mathbb E^4 }

\begin{document}

\title[Advances in Metallic Riemannian Geometry]{Recent Advances in Metallic Riemannian Geometry: A Comprehensive Review}

\author[B.-Y. Chen, M. A. Choudhary and A. Perween]{Bang-Yen Chen, Majid Ali Choudhary and Afshan Perween}

\begin{abstract} {Metallic structures, introduced by V. de Spinadel in 2002, opened a new avenue in differential geometry. Building upon this concept, C. E. Hre\c{t}canu and M. Crasmareanu laid the foundation for metallic Riemannian manifolds in 2013. The field's rich potential and diverse applications have since attracted significant research efforts, leading to a wealth of valuable insights. This review delves into the latest advances in metallic Riemannian geometry, a rapidly progressing area within the broader field of differential geometry.}
\end{abstract}

\keywords{Metallic Riemannian manifolds, metallic semi-Riemannian manifolds, metallic-like statistical manifolds, almost-complex metallic manifolds}

\subjclass[2020]{53C15, 53C25, 53C40, 53C42, 53C50, 11B39, 53B05, 53B20, 53B25, 53C20, 53C35, 53C75}

\maketitle

\section{Introduction}\label{S1}

In the realm of differential geometry, Riemannian manifolds have become a focal point of research in recent years. These curved spaces take centre stage when imbued with an intriguing structure inspired by the concept of metallic means. Spinadel \cite{2002} introduced this notion in 2002, generalizing the well-known golden ratio and encompassing a broader family of proportions like the silver, bronze, and copper mean. This extension allows mathematicians to explore the geometric properties of Riemannian manifolds with a richer and more nuanced perspective. Goldberg and Yano defined a broader category called polynomial structures \cite{yano}. Within this category, there's a specific type called a metallic structure. This metallic structure can be seen as a more general version of a special kind of polynomial structure, known as a golden structure. Within the field of differential geometry, the concept of a metallic structure applied to differentiable manifolds finds its origins in the work of Hre\c{t}canu and Crasmareanu \cite{3}. Their research extended the previously established theory of golden structures. Ozkan and Yılmaz in \cite{ozkan} explored metallic structures using almost product structures. Building on the concept of metallic Riemannian structures, Gezer and Karaman \cite{in} dive into their integrability conditions and curvature properties. Their analysis hinges on the application of a specifically chosen operator. Blaga and Hre\c{t}canu (2018) introduced the notion of metallic conjugate connections within the framework of metallic structures \cite{blaga2018metallic}. These connections represent a generalization of golden conjugate connections. Researchers have explored various forms of metallic structures (see \cite{2018, harmonic}).

Submanifolds occupy a particularly fascinating niche within the field of differential geometry. The fact that a submanifold of a Riemannian manifold is always a Riemannian one is widely recognized. The incorporation of a metallic structure on an ambient Riemannian manifold serves as a significant tool for deriving the geometric properties of its submanifolds. This approach has become a powerful technique in the investigation of submanifold geometry within Riemannian contexts. Employing a methodology similar to that used in the context of golden Riemannian manifolds, Hre\c{t}canu and Crasmareanu laid the groundwork for studying the differential geometry of metallic Riemannian manifolds in \cite{3}, particularly invariant submanifolds, which were characterized here. Blaga and Hre\c{t}canu in \cite{6} investigated a specific type of submanifold within metallic Riemannian manifolds called an invariant submanifold and proved that Nijenhuis tensor of the tensor field of the induced structure is identically zero on an invariant submanifold of locally decomposable metallic Riemannian manifolds, maintaining the property of the ambient manifold's locally decomposability. This submanifold was further extended to slant, semi-slant, hemi-slant, and bi-slant submanifolds, respectively, in metallic and golden Riemannian manifolds \cite{6, bl, hemi} by the same authors. Furthermore, in metallic Riemannian manifolds, the study of totally umbilical semi-invariant submanifolds was carried out in \cite{totally}.

The warped products can be thought of as a natural extension of cartesian products. Beginning with the research of Nash, this idea entered the field of mathematics. Nash established an embedding theorem, according to which any Riemannian manifold can be isometrically embedded into a Euclidean space. Moreover, Nash's theorem demonstrates that every warped product $\bar M_1 \times_f \bar M_2$ can be embedded in a Euclidean space as a Riemannian submanifold \cite{nash}. First, Bishop and Neill \cite{w2} created the idea of the warped product in Riemannian manifolds, which led to the creation of a broad class of complete manifolds with negative curvature. A surface of revolution gave rise to the idea of the warped product. From the definition of the warped product, N\"olker \cite{nolker} developed the idea of multiple warped products. A bi-warped product is a specific instance of a multiply-warped product. In addition to differential geometry, the warped product is a crucial concept in mathematical physics, specifically in general relativity. Warped products include the Robertson-Walker model, the Kruscal model, the Schwarzschild solution, and the static model. The warped products can be used to express a large number of accurate solutions to modified field equations and Einstein field equations.
Inspired by all this, the study of warped product submanifolds in metallic Riemannian manifolds is carried out in \cite{w1} by Blaga and Hre\c{t}canu. furthermore, the same authors have studied warped product pointwise bi-slant submanifolds
as well as warped product pointwise semi-slant or hemi-slant submanifolds within metallic Riemannian
manifolds in \cite{w3}. The study of warped product lightlike submanifolds induced in metallic semi-Riemannian manifolds is done by Shanker and Yadav in \cite{wl}.

On the other hand, it is widely understood in Riemannian geometry that the metric induced on a submanifold of a Riemannian manifold remains Riemannian. However, the induced metric of the semi-Riemannian metric of the ambient manifold is not consistently non-degenerate in the semi-Riemannian situation. The preceding instance yields a fascinating class of submanifolds known as lightlike submanifolds. The methods used to study the geometry of submanifolds in the Riemannian situation are inapplicable in the semi-Riemannian case due to the degeneracy of the induced metric on lightlike submanifolds. As a result, the classical theory fails when attempting to define any induced object on a lightlike submanifold. The primary challenge stems from the nonzero intersection between the tangent bundle and the normal bundle of a lightlike submanifold. To overcome the challenges encountered when researching lightlike submanifolds, Duggal and Bejancu \cite{duggal} devised a non-degenerate distribution known as the screen distribution to create a lightlike transversal vector bundle that does not intersect with its lightlike tangent bundle. It is commonly recognized that numerous significant characterizations in lightlike geometry can be obtained by selecting an appropriate screen distribution. This discovery sparked a wave of research, with mathematicians delving into various lightlike submanifolds within diverse manifold types. Their efforts yielded significant geometric insights and classifications. In 2018 Acet initiated the study of lightlike hypersurfaces in metallic semi-Riemannian manifolds \cite{acet}. The notion of lightlike submanifolds in metallic semi-Riemannian manifolds is introduced in \cite{li}. Later different types of lightlike submanifolds induced in metallic semi-Riemannian manifolds are studied (see \cite{tr, half, lone} etc).

 Furthermore, a central challenge in submanifold geometry is to establish an optimal inequality between the intrinsic and extrinsic invariants of a submanifold. In 1993, the first author introduced a groundbreaking concept called Chen invariants (or $\delta$-invariants) \cite{16}. These invariants provided a powerful tool for analyzing the relationship between a submanifold's internal properties (intrinsic) and its interaction with the larger space (extrinsic). Chen's work achieved this by establishing an optimal inequality. This breakthrough paved the way for a whole new area of differential geometry, sparking extensive research into Chen invariants and related inequalities for various submanifolds in diverse ambient spaces. Inspired by this, the second author with Uddin in \cite{17} has explored Chen-type Inequalities for slant submanifolds in metallic Riemannian space forms. The Chen-Ricci inequality for isotropic submanifolds in locally metallic product space forms was proven by Li et al. in \cite{29}. Moreover, the Chen-type inequality in metallic-like statistical manifolds was examined by Bahadir in \cite{B21}.

In submanifold geometry, a paradigm shift occurred with the replacement of the classical Gauss curvature by the Casorati curvature. This substitution provided the foundation for the establishment of optimal inequalities for submanifolds in diverse ambient spaces, using the power of Casorati curvatures \cite{11}. Casorati favoured this curvature above the conventional Gauss curvature as it more closely matches the perception of curvature because both principal curvatures of a surface in $\mathbb E^3$ are zero if and only if the Casorati curvature vanishes. Numerous scholars have employed the Casorati curvature technique extensively to generate optimal inequalities for submanifolds in metallic Riemannian manifolds. In \cite{14}, Choudhary and Blaga examined sharp inequalities using generalized normalized $\delta$-Casorati curvatures in a metallic Riemannian space.

 Wintgen is credited with discovering the Wintgen inequality \cite{18}. This inequality is a significant result in differential geometry because it relates the inherent properties (intrinsic invariants) of a surface $M^2$ in 4-dimensional Euclidean space $\mathbb E^2$ to extrinsic invariants. 

Given recent advances in metallic differential geometry, this paper provides an extensive overview of the most important recent developments in the subject. Explores the most recent studies, examining how these results shed light on the complex connection between metallic structures and multiple geometric features. With the help of this critical examination, researchers should be able to grasp the fundamental concepts guiding the fascinating new developments in metallic differential geometry.

\section{Preliminaries}\label{S2}

\subsection{Metallic Riemannian manifolds}\label{S2.1}

 Consider $(\bar M, g)$ as an $m$-dimensional Riemannian manifold. A $(1,1)$-tensor field $F$ on $\bar M$ is called a polynomial structure \cite{yano} if it satisfies the condition $B(F) = 0$, where $$B(X):= X^n +a_n X^{n-1} + \cdots +a_2 X+a_1 I,$$ with any real numbers $a_1,a_2,\cdots, a_n,$ where $I$ is the identity transformation on $\Gamma (T \bar M )$. One can observe that 
\begin{itemize}
    \item $F$ is an almost complex structure when $B(X) = X^2 + I$,
    \item $F$ is an almost product structure when $B(X) = X^2 - I$.
\end{itemize}

Again, suppose that $(\bar M,g)$ is an $m$-dimensional Riemannian manifold. Then for any integer $p,q$ the $(1,1)$-tensor field $\varphi$ on $\bar M$ is called a {\it metallic structure} \cite{3} if $$\varphi ^2 =p \varphi +qI,$$ where $I$ denotes the identity operator on $\Gamma (T \bar M )$. It is commonly recognized that a Riemannian metric $g$ is termed $\varphi$-compatible if the following relation is satisfied $$g(X, \varphi Y)=g(\varphi X, Y)$$ for any $X, Y \in \Gamma(T \bar M).$
A manifold $\bar M$ that is equipped with a metallic structure $\varphi$ and a $\varphi$-compatible Riemannian metric $g$ is called a {\it metallic Riemannian manifold}.
Now, replacing $X$ by $\varphi X,$ we have $$g(\varphi X , \varphi Y)=pg(X, \varphi Y)+qg (X, Y).$$
Specifically, when $ p = q = 1,$ a metallic Riemannian manifold is referred to as a {\it golden Riemannian manifold} \cite{MH08,HC09}.

A $(1, 1)$-tensor field $F$ is referred to as an almost product structure \cite{4} on a Riemannian manifold $(\bar M, g)$ of dimension $m$ if it meets the conditions $F^2 = I$ and $F \neq \pm I$. Furthermore, $(\bar M, g)$ is described as an almost product Riemannian manifold when the almost product structure $F$ satisfies $ \forall X, Y \in \Gamma(T \bar M)$ $$g(FX, Y ) = g(X, F Y ).$$ 

Via any metallic structure $\varphi$ on $\bar M$ one obtains two almost product structures on $\bar M$ \cite{3}
\begin{equation} \label{1}
  \begin{aligned}
    F_1 &= \frac{2}{2\sigma_{p,q} - p} \varphi - \frac{p}{2\sigma_{p,q} - p} I, \\
    F_2 &= -\frac{2}{2\sigma_{p,q} - p} \varphi + \frac{p}{2\sigma_{p,q} - p} I.
  \end{aligned}
\end{equation}
Here $\sigma_{p,q} = \frac{p+ \sqrt{p^2 +4q}}{2}$ symbolize the metallic proportions or the members of the metallic means family. Additionally, any almost product structure $F$ on $\bar M$ yields two metallic structures 
\begin{equation}
\begin{aligned}
    \varphi_1 &= \frac{p}{2} I +\frac{2\sigma_{p,q} -p}{2} F, \\
    \varphi_2 &= \frac{p}{2} I -\frac{2\sigma_{p,q} -p}{2} F.
\end{aligned}
\end{equation}
 
\begin{definition} 

We recall the following definition from \cite{6}.

{\rm (i)} Assume that $(\bar M, g, \varphi)$ denotes a metallic Riemannian manifold and let $\nabla$ be a linear connection on $\bar M$. The linear connection $\nabla$ is referred to as a $\varphi$-connection if $\varphi$ remains covariantly constant with respect to $\nabla$, meaning it satisfies $\nabla \varphi = 0.$

{\rm (ii)} A metallic Riemannian manifold \((\bar{M}, g, \varphi)\) is described as a locally metallic Riemannian manifold if the Levi-Civita connection \(\bar{\nabla}\) of \(g\) is a \(\varphi\)-connection.
\end{definition} 

\begin{example}
With the canonical coordinates $\left(\zeta_1, \ldots, \zeta_7\right)$ and the natural Euclidean metric $\langle\cdot, \cdot\rangle$, let us consider the Euclidean $7$-space $\mathbb{E}^7$. We then define the immersion $f: \bar{M} \rightarrow \mathbb{E}^7$ by

$$
f\left(t_1, t_2\right):=\left(\sin t_1+2, \cos t_1, 2 \sin t_2+1, 2 \cos t_2, t_2, 2 t_1, 3\right),
$$
where $\bar{M}:=\left\{\left(t_1, t_2\right) \mid t_1, t_2 \in\left(0, \frac{\pi}{2}\right)\right\}$.
Then a local orthonormal frame on $T \bar{M}$ is given by the vector fields
$$
\begin{aligned}
& \zeta_1=\cos t_1 \frac{\partial}{\partial \zeta_1}-\sin t_1 \frac{\partial}{\partial \zeta_2}+2 \frac{\partial}{\partial \zeta_6}, \\
& \zeta_2=2 \cos t_2 \frac{\partial}{\partial \zeta_3}-2 \sin t_2 \frac{\partial}{\partial \zeta_4}+\frac{\partial}{\partial \zeta_5} .
\end{aligned}
$$
We consider the metallic structure $\varphi:\mathbb{E}^7 \rightarrow \mathbb{E}^7$ :
$$
\varphi\left(\frac{\partial}{\partial \zeta_i}\right)=\sigma \frac{\partial}{\partial \zeta_i}, \quad i \in\{1,2,5\} ; \quad \varphi\left(\frac{\partial}{\partial \zeta_j}\right)=\bar{\sigma} \frac{\partial}{\partial \zeta_j}, j \in\{3,4,6,7\},
$$
where $\sigma:=\sigma_{p, q}=\frac{p+\sqrt{p^2+4 q}}{2}$ is the metallic number, $p, q$ are integers, and $\bar{\sigma}=p-\sigma$. Then $\left(\mathbb{E}^7,\langle\cdot, \cdot\rangle, \varphi\right)$ is a metallic Riemannian manifold. 
\end{example}

Let $\left(\bar M=M_p\left(c_p\right) \times M_q\left(c_q\right), F\right)$ be a locally Riemannian product manifold, where $M_p$ and $M_q$ have constant sectional curvatures $c_p$ and $c_q$, respectively. Then, the Riemannian curvature tensor $R$ of $\bar M=M_p\left(c_p\right) \times M_q\left(c_q\right)$ for $X,Y$ and $Z \in \Gamma (T\bar M)$ is \cite{7}
\begin{equation} \label{3}
\begin{aligned} 
R(X, Y) Z  =&\, \frac{1}{4} \left( c_p + c_q \right) [g(Y, Z) X - g(X, Z) Y + g(FY, Z) FX \\
& - g(FX, Z) FY] + \frac{1}{4} \left( c_p - c_q \right) [g(FY, Z) X \\
& - g(FX, Z) Y + g(Y, Z) FX - g(X, Z) FY].
\end{aligned}
\end{equation}

In view of \eqref{1} and \eqref{3}, one can get \cite{8}
\begin{equation}
\begin{aligned}
R(X, Y) Z = &\, \frac{1}{4}\left(c_p+c_q\right)[g(Y, Z) X-g(X, Z) Y] \\
& +\frac{1}{4}\left(c_p+c_q\right)\left\{\frac{4}{\left(2 \sigma_{p, q}-p\right)^2}[g(\varphi Y, Z) \varphi X-g(\varphi X, Z) \varphi Y]\right. \\
& +\frac{p^2}{\left(2 \sigma_{p, q}-p\right)^2}[g(Y, Z) X-g(X, Z) Y] \\
& +\frac{2 p}{\left(2 \sigma_{p, q}-p\right)^2}[g(\varphi X, Z) Y+g(X, Z) \varphi Y \\
& -g(\varphi Y, Z) X-g(Y, Z) \varphi X]\Bigg\} \\
& \pm \frac{1}{2}\left(c_1-c_2\right)\left\{\frac{1}{2 \sigma_{p, q}-p}[g(Y, Z) \varphi X-g(X, Z) \varphi Y]\right. \\
& +\frac{1}{2 \sigma_{p, q}-p}[g(\varphi Y, Z) X-g(\varphi X, Z) Y] \\
& \left.+\frac{p}{2 \sigma_{p, q}-p}[g(X, Z) Y-g(Y, Z) X]\right\}. 
\end{aligned}
\end{equation}

\subsection{Metallic semi-Riemannian manifolds}\label{S2.2}

The positive solution of $$x^2 -px-q=0,$$ is called a member of the metallic means family \cite{2002}, where $p, q$ are fixed positive integers. These numbers denoted by; $$\sigma_{p, q}=\frac{p+\sqrt{p^2+4 q}}{2,}$$ are known {\it $(p, q)$-metallic numbers}.

A polynomial structure on a semi-Riemannian manifold $\bar M$ is called metallic if it is determined by $\varphi$ such that $$\varphi^2 = p\varphi + qI.$$
If a semi-Riemannian metric $g$ satisfies the equation
$$g(X,\varphi Y) =g(\varphi X,Y),$$
which yields
$$g(\varphi X,\varphi Y) = pg(X,\varphi Y)+ qg(X,Y),$$
then $g$ is called $\varphi$-compatible.

\begin{definition} \cite{acet}
A semi-Riemannian manifold $(\bar M,g)$ equipped with $\varphi$ such that the semi-Riemannian metric $g$ is $\varphi$-compatible is named a metallic semi-Riemannian manifold and $(g,\varphi)$ is called a metallic semi-Riemannian structure on $\bar M$.
\end{definition}

\subsection{Metallic-like statistical manifolds}\label{S2.3}

Takano proposed generalized almost complex and almost contact statistical manifolds, referring to them as K\"ahler-like statistical manifolds and Sasaki-like statistical \cite{9}. Motivated by this research, Bahadır in \cite{B21} elucidated metallic-like statistical manifolds, a generalized form of metallic manifolds.

\begin{definition} \cite{B21}
    Consider a locally metallic semi-Riemannian manifold $(\bar M, g, \varphi)$ equipped with a tensor field $\varphi^*$ of type $(1,1)$, satisfying $$g(\varphi X, Y) = g(X, \varphi^* Y)$$ for any vector fields $X,Y$. From the above equation one can easily derive $$(\varphi)^* X = p \varphi^* X+qX,$$ $$g(\varphi X, \varphi^* Y)=pg(\varphi X,Y)+qg(X,Y).$$ Then $(\bar M, g, \varphi)$ is called {\it metallic-like statistical manifold}.

    According to the above two equations, the tensor fields $\varphi +\varphi^*$ and $\varphi -\varphi^*$ are symmetric and skew-symmetric concerning $g$, respectively. 
\end{definition}

\begin{proposition} {\rm\cite{B21}}
$(\bar M, g, \varphi)$ is a metallic-like statistical manifold if and only if so is $(\bar M, g, \varphi^*).$
\end{proposition}

If one chooses $\varphi=\varphi^*$ in a metallic-like statistical manifold, then we have a metallic semi-Riemannian manifold.

\begin{example} \cite{B21}
Consider the semi-Euclidean space $\mathbb E^3_1$ with standard coordinates $(\zeta_1 ,\zeta_2 ,\zeta_3)$ and the semi-Riemannian
metric $g$ with the signature $(-, +, +)$. Let $\varphi$ be a $(1,1)$-tensor field on $\mathbb E^3_1$ defined by
\[
\varphi(\zeta_1, \zeta_2, \zeta_3) = \frac{1}{2}\big(p\zeta_1 + (2\sigma - p)\zeta_2, p\zeta_2 + (2\sigma - p)\zeta_1, (p-\sigma)\zeta_3\big)
\]
for any vector field $(\zeta_1 , \zeta_2 , \zeta_3) \in \mathbb E^3_1$ where $\sigma:=\sigma_{p, q}=\frac{p+\sqrt{p^2+4 q}}{2}$ are members of the metallic means family. Then we have $\varphi^2 =p \varphi +qI$. Also, we can easily see that the structure is compatible with the metric. This implies that $\varphi$ is a metallic structure on $\mathbb E^3_1$.
Now define a $(1,1)$-tensor field $\varphi^*$ on $\mathbb E^3_1$ by
\[
\varphi^* (\zeta_1, \zeta_2, \zeta_3) = \frac{1}{2} \left(p\zeta_1 + (p-2\sigma )\zeta_2, p\zeta_2 + (p-2\sigma)\zeta_1, (p-\sigma)\zeta_3\right).
\]
Thus, we find $\varphi^{*2} =p \varphi^* +qI$. Furthermore, we also have $g(\varphi X, Y) = g(X, \varphi^* Y)$. Consequently, $(\mathbb E^3_1,g,\varphi)$ is a metallic-like statistical manifold.
\end{example}

\subsection{Almost-complex metallic manifolds}\label{S2.4}

Consider the following: $$x^2 -ax+\frac{3}{2}b,$$ where $a$ and $b$ are the real numbers that satisfy $-\sqrt{6b}<a<\sqrt{6b}$ and $b \geq 0$. The positive solution of the equation has complex roots given as:
\begin{equation} \label{acm}
    \mathcal{C}_{a,b} =\frac{a+\sqrt{a^2 -6b}}{2},
\end{equation}
which is named as complex metallic means family in \cite{acm}. In particular, if $a = 1$ and $b = 1,$ then the complex metallic means family $\mathcal{C}_{a,b} =\frac{a+\sqrt{a^2 -6b}}{2}$
reduces to the complex golden mean: $\mathcal{C}_{1,1}=\frac{1+\sqrt{5i}}{2},i^2=-1.$

In \cite{acm}, using the complex mean given in \eqref{acm}, the authors defined a new type of structure on a Riemannian manifold. Let $\bar M$ be a Riemannian manifold. An almost complex metallic structure is a $(1, 1)$-tensor field, $\varphi_{\bar M}$ on $\bar M$, which satisfies the relation
\begin{equation} \label{acm2}
    \varphi^2_{\bar M}-a\varphi_{\bar M}+\frac{3}{2}bI=0,
\end{equation}
where $I$ is the identity $(1, 1)$-tensor field on $\bar M$. In this case, $\bar M$ is referred to as an almost complex metallic manifold equipped with an almost complex structure $\varphi_{\bar M}$.

Note that if we take $a = m$ and $b = \frac{2}{3}$
in \eqref{acm2}, we obtain an almost poly-Norden structure.

\begin{example} \cite{acm2}
Cosider the $4$-tuples real space $\mathbb R^4$. Let $\varphi_{\bar M}: \mathbb R^4 \to \mathbb R^4$ 
be a map given by $$\varphi_{\bar M} (\zeta_1 ,\zeta_2 ,\omega_1 ,\omega_2)=\big(\mathcal{C}_{a,b}\zeta_1 ,\mathcal{C}_{a,b}\zeta_2 ,(a-\mathcal{C}_{a,b})\omega_1 ,(a-\mathcal{C}_{a,b})\omega_2\big),$$ where $\mathcal{C}_{a,b} =\frac{a+\sqrt{a^2 -6b}}{2}.$ 
One can easily see that $\varphi_{\bar M}$ satisfies $\varphi^2_{\bar M}-a\varphi_{\bar M}+\frac{3}{2}bI=0$. Thus, $(\mathbb R^4 ,\varphi_{\bar M})$ is an example of almost complex metallic manifolds.

If $(\bar M,g)$ is a semi-Riemannian manifold equipped with an almost complex metallic structure such that the metric $g$ is $\varphi_{\bar M}$-compatible, $$g(\varphi_{\bar M} X,\varphi_{\bar M} Y)=ag(\varphi_{\bar M} X, Y)-\frac{3}{2}bg(X, Y)$$ equivalent to $$g(\varphi_{\bar M} X, Y)=g(X,\varphi_{\bar M} Y)$$ for every $X, Y \in \Gamma(T\bar M).$ Hence, $(\bar M,\varphi_{\bar M},g)$  is called an almost-complex metallic semi-Riemannian (briefly ACMSR) manifold.
\end{example}

\begin{remark}
{\rm Inspired by the Meta-Golden-Chi ratio, the authors in \cite{meta} developed the concept of Meta-Metallic manifolds via the Meta-Metallic-Chi ratio and Metallic manifolds. They also provided an example, explored the specific characteristics of the Meta-Metallic structure, and derived the conditions necessary for the integrability of the almost Meta-Metallic structure, as well as obtained its relationship with the curvature tensor.}
\end{remark}

\section{Submanifolds Immersed of Metallic Riemannian \\Manifolds}\label{S3}

\subsection{Invariant submanifolds in metallic Riemannian manifolds}\label{S3.1}

A submanifold $M$ of $\bar{M}$ is called invariant if $\varphi(T_x M) \subset T_x M$ for all $x \in M$. Hence we get $\varphi(T_x M^\perp) \subset T_x M^\perp$ for any $x \in M$, since we have $g(X, \varphi U) = g(\varphi X, U) = 0$ for any $X \in \Gamma(TM)$ and any $U \in \Gamma(TM^\perp)$ \cite{6}.

Many researchers have studied invariant submanifolds in metallic Riemannian manifolds. In \cite{6}, Blaga and Hre\c{t}canu studied the properties of invariant isometrically immersed submanifolds in these manifolds, with a particular emphasis on the induced $\Sigma$-structure.

\begin{proposition}\label{P:3.1} {\rm \cite{6}}
Given a locally metallic Riemannian manifold $(\bar M, \varphi, g),$ with $M$ an isometrically immersed invariant submanifold, then for every $X, Y \in \Gamma(TM)$ we have: \begin{align*}&\nabla \varphi =0,\\&\sum_{\alpha=1}^{r} h_{\alpha}(X, \varphi Y)\xi_{\alpha} = \sum_{\alpha=1}^{r} h_{\alpha}(X, Y)\varphi \xi_{\alpha} = \sum_{\alpha=1}^{r} h_{\alpha}(\varphi X, Y)\xi_{\alpha},\\ &h_{\alpha}(\varphi X, \varphi Y ) = ph_{\alpha} (X, \varphi Y ) + qh_{\alpha} (X, Y ) \;\; \text{for any} \;\;  1 \leq \alpha \leq r,\end{align*}
where the (symmetric) second fundamental tensors corresponding to $\xi_{\alpha}$ are $h_{\alpha}, 1 \leq \alpha \leq r,$ and $\xi_{\alpha}$ is an orthonormal basis for the normal space. i.e., $h(X, Y ) = \sum_{\alpha=1}^{r} h_{\alpha}(X, Y )\xi_{\alpha},$ for $X, Y \in \Gamma(TM).$
\end{proposition}

\begin{proposition} \label{P:3.2}{\rm \cite{6}}
A locally metallic Riemannian manifold $(\bar M, \varphi, g)$ of dimension $n+r$ has an isometrically immersed invariant n-dimensional submanifold of codimension r, denoted by $M$ and $\mathbf{\Sigma}:= (T, g, \eta_{\alpha} = 0, \xi'_{\alpha} = 0, (a_{\alpha\beta}))_{1 \leq \alpha, \beta \leq r}$ is the induced structure in $M$. Where $T$ is the tangent bundle in $M,$ $\xi'_\alpha$ is a vector field on $M$,$\eta_{\alpha}$ is a 1-form in $M$ and $(a_{\alpha\beta}))_{1 \leq \alpha, \beta \leq r}$ is a $r \times r$ matrix of smooth real function on $M$. Then $$TS_{\xi_\alpha} = S_{\xi_\alpha} T,$$
the Nijenhuis tensor field of $T$ disappears identically in $M$ for any $1 \leq \alpha \leq r$.(i.e.,
$\mathcal{N}_T (X, Y ) = 0,$ for any $X, Y \in \Gamma(TM))$.
\end{proposition}

\begin{proposition} \label{P:3.3} {\rm \cite{6}}
Consider $M$ as an isometrically immersed invariant n-dimensional submanifold of codimension r of the locally metallic $(n + r)$-dimensional Riemannian manifold $(\bar M, \varphi, g)$ and $\mathbf{\Sigma}:= (T, g, \eta_{\alpha}= 0, \xi'_{\alpha}=0, (a_{\alpha\beta}))_{1 \leq \alpha, \beta \leq r}$ is the induced structure on $M$. 
Then the components $\mathcal{N}^{(2)} , \mathcal{N}^{(3)} and \mathcal{N}^{(4)}$ vanish identically on $M$. Moreover, if $\mathcal{N}_T = 0,$ then $\mathcal{N}^{(1)}$ vanishes, too, on $M$. This specifically occurs if the normal connection $\nabla^\perp$ on the normal bundle vanishes in identical way i.e., $ \lambda_{\alpha\beta} = 0,$ for every $1 \leq \alpha, \beta \leq r,$ where $\lambda_{\alpha\beta} =-\lambda_{\beta\alpha}$ is $1$-form on $M$ corresponding to the normal connection $\nabla^\perp$.
\end{proposition}

\begin{remark}
  {\rm In \cite{mg}, Gök has examined invariant submanifolds within metallic Riemannian manifolds and has derived results indicating that an isometrically immersed submanifold of a metallic Riemannian manifold is a totally geodesic invariant submanifold.} 
\end{remark}
  
  Gök and Kılıç have investigated totally umbilical semi-invariant submanifolds in locally decomposable metallic Riemannian manifolds, as reported in \cite{totally}, presenting the following conclusions.

\begin{proposition} \label{P:3.4} {\rm \cite{totally}}
Let $M$ be a totally umbilical proper semi-invariant submanifold within a locally decomposable metallic Riemannian manifold $(\bar M,g,\varphi )$. The following statements hold:

{\rm (i)} The invariant distribution $D^\theta$ is integrable.

{\rm (ii)} The anti-invariant distribution $D^\perp$ is integrable.
\end{proposition}

\begin{proposition} \label{P:3.5} {\rm \cite{totally}}
For any totally umbilical proper semi-invariant submanifold $M$ of a locally decomposable metallic Riemannian manifold $(\bar M,g,\varphi )$, the covariant derivative of the endomorphism $T$ is zero.
\end{proposition}

\begin{theorem} \label{T:3.1} {\rm \cite{totally}}
Consider $M$ as a totally umbilical proper semi-invariant submanifold of a locally decomposable metallic Riemannian manifold $(\bar M,g,\varphi )$. If $\dim \bar M = \dim M + \dim D^\perp$, then $M$ is a totally geodesic submanifold.
\end{theorem}

\begin{theorem}\label{T:3.2} {\rm \cite{totally}}
If $M$ is a totally umbilical proper semi-invariant submanifold with a non-zero mean curvature vector $\mathcal{H}$ in a locally decomposable metallic Riemannian manifold $(\bar M,g,\varphi)$, and if $\dim D^\theta \geq 2$, then $M$ is an extrinsic sphere.
\end{theorem}

\begin{theorem}\label{T:3.3} {\rm \cite{totally}}
In a locally decomposable metallic Riemannian manifold $(\bar M,g,\varphi)$ with positive or negative curvature, there are no totally umbilical proper semi-invariant submanifolds.
\end{theorem}

\subsection{Anti-invariant submanifolds of metallic Riemannian manifolds}\label{S3.2}

For any $x \in M$, a submanifold $M$ of $\bar M$ is said to be anti-invariant if $\varphi (T_xM) \subset T_x M^\perp$.

\begin{proposition} \label{P:3.6} {\rm \cite{6}}
Given a locally metallic Riemannian manifold $(\bar M,\varphi,g)$ and an isometrically immersed anti-invariant submanifold $M$, then for every $X, Y \in \Gamma(TM)$
$$\sum_{\alpha=1}^{r} h_{\alpha}(X, Y)t \xi_{\alpha} = -\sum_{\alpha=1}^{r}g(\varphi Y, \xi_{\alpha}) S_{\xi_\alpha}X,$$ $$\sum_{\alpha=1}^{r} h_{\alpha}(X, Y)n \xi_{\alpha} =\sum_{\alpha=1}^{r}g(\varphi Y, \xi_{\alpha}) \nabla^\perp_X \xi_{\alpha} +\sum_{\alpha=1}^{r} X(g(\varphi Y,\xi_{\alpha}))-\varphi (\nabla_X Y),$$
where $t : \Gamma(TM^\perp ) \rightarrow \Gamma(TM),\, tU = (\varphi U)^T$ and $n : \Gamma(TM^\perp) \rightarrow \Gamma(TM^\perp),\, nU = (\varphi U)^\perp .$
\end{proposition}

\begin{proposition} \label{P:3.7} {\rm \cite{6}}
For a (n+r)-dimensional locally metallic Riemannian manifold $(\bar M, \varphi, g)$, let $M$ be an isometrically immersed anti-invariant n-dimensional submanifold with codimension r and $\mathbf{\Sigma}:= (T, g, \eta_{\alpha}= 0, \xi'_{\alpha}=0, (a_{\alpha\beta}))_{1 \leq \alpha, \beta \leq r}$ is the induced structure on $M$.  Then the components $\mathcal{N}^{(2)}$ and $\mathcal{N}^{(3)}$ vanish identically on $M$. Furthermore, if $\xi'_\alpha$ are parallel for a linear symmetric connection, for any $1 \leq \alpha \leq r,$ then $\mathcal{N}^{(1)}$ and $\mathcal{N}^{(4)}$ also vanish on $M$. 
\end{proposition}

\subsection{Slant submanifolds of metallic Riemannian manifolds}\label{S3.3}

A submanifold $M$ is called {\it slant} if the angle $\theta (X_x)$ between $\varphi (X_x)$ and $T_x M$ is independent of the choice of the point $x \in M$ and the non-zero tangent vector $X$ at $x$. Here, $\theta =: \theta (X_x)$ is called the slant angle. Slant submanifolds include invariant and anti-invariant submanifolds, where the slant angles are $\theta = 0$ and $\theta = \frac{\pi}{2}$, respectively. A slant immersion that is neither invariant nor anti-invariant is called a proper slant  \cite{6}.

Blaga and Hre\c{t}canu  \cite{6} investigated the characteristics of slant isometrically immersed submanifolds
in metallic manifolds, particularly focusing on the induced $\Sigma$-structure.

\begin{proposition} \label{P:3.8} {\rm \cite{6}}
Let $M$ be an isometrically immersed submanifold of the metallic Riemannian manifold $(\bar M, \varphi, g)$. If $M$ is slant with the slant angle $\theta$, then we have: $$g(T X, T Y ) = \cos^2 {\theta}[pg(X, \varphi Y ) + qg(X, Y )]$$ and
$$g(NX, NY ) = \sin^2 {\theta}[pg(X, \varphi Y ) + qg(X, Y )]$$
for any $X, Y \in \Gamma(TM)$. Where $TX \in \Gamma(TM)$ and $NX \in \Gamma(TM^\perp)$ are tangential and normal parts of $\varphi X , \varphi V$ respectively, for any $V \in \Gamma(TM^\perp).$
\end{proposition}

\begin{proposition} \label{P:3.9} {\rm  \cite{6}}
Let $M$ be an isometrically immersed submanifold of the metallic Riemannian manifold $(\bar M, \varphi, g)$. If $M$ is slant with the slant angle $\theta$, then: $$(\nabla_X T^2 )Y = p \cos^2 {\theta}(\nabla_X T)Y$$
for any $X, Y \in \Gamma(TM).$   
\end{proposition}

\begin{corollary} \label{C:3.4} {\rm \cite{6}}
If $M$ is an isometrically immersed slant submanifold of the metallic
Riemannian manifold $(\bar M, \varphi, g)$ with the slant angle $\theta$, then $\nabla^2 T=0$ if and only if $M$ is anti-invariant or $(M, T, g)$ is locally metallic Riemannian manifold.
\end{corollary}

Hre\c{t}canu and Blaga have also introduced the following theorem in \cite{bl} for slant submanifold in metallic Riemannian manifold. 

\begin{theorem} \label{T:3.5} 
Let $M$ be a submanifold in the Riemannian
manifold $(\bar M,g)$ endowed with an almost product structure $F$
on $\bar M$ and let $\varphi$ be the induced metallic structure by $F$ on $(\bar M,g)$. If $M$ is a slant submanifold in the almost product Riemannian manifold  $(\bar M,g,F)$ with the slant angle $
\theta$ and $F \neq -I$ $(I$ is the
identity on $\Gamma(TM))$ and $\varphi =\frac{((2\sigma_{p,q}-p)}{2)}F+(\frac{p}{2}I),$ then $M$ is a slant submanifold in the metallic Riemannian manifold
$(\bar M,g,\varphi )$ with slant angle $\theta$ given by $$\sin{\Theta} = \frac{2\sigma_{p,q}-p}{2\sigma_{p,q}}\sin{\theta}.$$    
\end{theorem}

\subsection{Bi-slant submanifolds of metallic Riemannian manifolds}\label{S3.4}

Let $M$ be an immersed submanifold in a metallic Riemannian manifold $(\bar M,g,\varphi)$. Then $M$ is called {\it bi-slant} if there exist two orthogonal differentiable distributions $D^\theta$ and $D^\perp$ on $M$ such that $TM = D^\theta \oplus D^\perp$ and $D^\theta, D^\perp$ are slant distributions with the slant angles $\theta_1$ and $\theta_2,$ respectively \cite{bl}.

The following result is obtained by the authors of \cite{bl} for bi-slant submanifolds in metallic Riemannian manifolds.

\begin{proposition} \label{P:3.10} 
If $M$ is a bi-slant submanifold in a metallic
Riemannian manifold $(\bar M,g,\varphi)$, with the slant angles $\theta_1 =\theta_2 =\theta$ and $g(\varphi X,Y) =0$ for any $X \in \Gamma(D^\theta)$
and $Y \in \Gamma (D^\perp),$ then $M$ is a slant submanifold in the metallic Riemannian manifold $(\bar M,g,\varphi)$ with the slant angle $\theta$.
\end{proposition}

\subsection{Semi-slant submanifolds of metallic Riemannian manifolds}\label{S3.5}

An immersed submanifold $M$ in a metallic Riemannian manifold $(\bar M,g,\varphi )$ is termed {\it semi-slant} if  the orthogonal distributions $D^\theta$ and  $D^\perp$ on $M$ satisfy the following conditions:

{\rm (i)} $TM=D^\theta \oplus D^\perp ;$

{\rm (ii)} $\varphi D^\theta =D^\theta;$

{\rm (iii)} $D^\perp$ is slant with $\theta \neq 0.$

\noindent Moreover,  if $\dim(D^\theta)  \dim(D^\perp) \neq 0$, then $M$ is a proper semi-slant submanifold \cite{bl}.

Hre\c{t}canu and Blaga in \cite{bl} have discussed some properties of semi-slant submanifolds in metallic Riemannian manifolds.

\begin{proposition} \label{P:3.11} {\rm \cite{bl}}
If $M$ is a semi-slant submanifold of the metallic Riemannian manifold $(\bar{M}, g, \varphi)$ with slant angle $\theta$ corresponding to the distribution $D^\perp$, then we obtain
$$g(TP^\perp X,TP^\perp Y)=\cos^2 {\theta}[pg(TP^\perp X,P^\perp Y)+qg(P^\perp X,P^\perp Y)],$$  $$g(NX,NY)=\sin^2{\theta}[pg(TP^\perp X,P^\perp Y)+qg(P^\perp X,P^\perp Y)]$$ for any $X,Y \in \Gamma(TM),$ where $P^\perp$ is the orthogonal projection on $D^\perp$ and $TX \in \Gamma(TM)$ and $NX \in \Gamma(TM^\perp)$ are the tangential and normal parts of $\varphi X , \varphi V$ respectively, for any $V \in \Gamma(TM^\perp).$
\end{proposition}

\begin{proposition} \label{P:3.12} {\rm \cite{bl}}
Let $M$ be a semi-slant submanifold of the
metallic Riemannian manifold $(\bar M,g,\varphi)$ with the slant angle $\theta$ of the distribution $D^\perp$. Then $$(TP^\perp)^2 =\cos^2{\theta}(pTP^\perp+qI),$$ where $I$ is the identity on $\Gamma(D^\perp)$ and $$\nabla((TP^\perp)^2)=p\cos^2{\theta} \nabla(TP^\perp).$$ 
\end{proposition}

\begin{remark}
     {\rm Examples of semi-slant submanifolds of the metallic manifold are provided in \cite{bl} and the authors have obtained integrability conditions for the distributions involved in the semi-slant submanifolds of Riemannian manifolds endowed with a metallic structure.}
\end{remark}

\subsection{Hemi-slant submanifolds of metallic Riemannian manifolds}\label{S3.6}

An immersed submanifold $M$ within metallic Riemannian manifold $(\bar M, g, \varphi)$ is termed as {\it hemi-slant} if orthogonal distributions $D^\theta$ and $D^\perp$ on $ M$ satisfy the following conditions:

{\rm (i)} $TM = D^\theta \oplus D^\perp ;$

{\rm (ii)} $D^\theta$ is slant with $\theta \in [0,\frac{\pi}{2}];$

{\rm (iii)} $D^\perp$ is anti-invariant distribution, i.e. $\varphi (D^\perp) \subseteq \Gamma(TM^\perp)$.

\noindent Moreover, if $\dim(D^\theta) \dim(D^\perp) \neq 0$ and $\theta \in (0, \frac{\pi}{2}),$ then $M$ is a proper hemi-slant submanifold \cite{hemi}.

\begin{theorem} \label{T:3.6} {\rm \cite{hemi}}
If $M$ is a hemi-slant submanifold within metallic Riemannian manifold $(\bar M, g, \varphi)$ then for any $X, Y \in \Gamma(TM)$ we have:
$$g(TP^\theta X,TP^\theta Y) =\cos^2{\theta}\, \{pg(TP^\theta X,P^\theta Y)+qg(P^\theta X,P^\theta Y)\},$$  $$g(NX,NY)=\sin^2{\theta}\,\{pg(TP^\theta X,P^\theta Y)+qg(P^\theta X,P^\theta Y)\},$$
where $P^\theta$ is the orthogonal projection on $D^\theta$, $TX \in \Gamma(TM)$ and $NX \in \Gamma(TM^\perp)$ are the tangential and normal parts of $\varphi X , \varphi V$ respectively, for any $V \in \Gamma(TM^\perp).$
\end{theorem}

\begin{theorem} \label{T:3.7} {\rm \cite{hemi}}
Consider a hemi-slant submanifold $M$ in a metallic Riemannian manifold $(\bar M, g, \varphi)$ with slant angle $\theta$ of the distribution $D^\theta$. Then:
$$(TP^\theta)^2 = \cos^2{\theta}\, (pTP^\theta + qI),$$
where $I$ is the identity on $\Gamma(D^\theta)$ and $$\nabla((TP^\theta)^2)=p\cos^2{\theta} \,\nabla(TP^\theta).$$  
\end{theorem}

\begin{theorem} \label{T:3.8} {\rm \cite{hemi}}
Consider $M$ to be an immersed submanifold within metallic Riemannian manifold $(\bar M, g, \varphi)$. Then $M$ is a hemi-slant submanifold in $\bar M$ if and only if  there exists a constant $c \in [0, 1]$ such that $D = \{X \in \Gamma(TM)|T^2 X = c(pT X + qX)\}$ is a distribution and $T Y = 0$ for any $Y$ perpendicular  to $D, Y \in \Gamma(TM),$ where $p, q$ are integers.
\end{theorem}

\begin{remark}
    {\rm Examples of hemi-slant submanifolds in metallic Riemannian manifold, the conditions for the integrability of the distributions of a hemi-slant submanifold in a metallic  Riemannian manifold and conditions for these submanifolds to be mixed totally geodesic is discussed in \cite{hemi}.}
\end{remark}

The notion of quasi-hemi-slant submanifolds of metallic Riemannian manifolds is studied  by Karmakar and Bhattacharyya  \cite{quasi}.

\begin{definition} {\rm \cite{quasi}}
A {\it quasi hemi-slant submanifold} $M$ of a metallic Riemannian manifold $(\bar{M} , g, \varphi)$ is a submanifold that admits three orthogonal complementary distributions $D, D^\theta, D^\perp$ such that 

{\rm (i)} $TM=D \oplus D^\theta \oplus D^\perp$.

{\rm (ii)} $D$ is invariant, that is, $\varphi D = D$.

{\rm (iii)} $D^\theta$ is slant with $\theta$, and hence $\theta$ is called the slant angle.

{\rm (iv)} $D^\perp$ is anti-invariant, i.e., $\varphi D^\perp \subseteq T^\perp M$.

In the above case, $\theta$ is called the quasi hemi-slant angle of $M$, and $M$ is called proper if $D \neq 0$, $D^\theta \neq 0$, $D^\perp \neq 0$, and $\theta \neq (0,\frac{\pi}{2})$.
\end{definition}

\begin{theorem} \label{T:3.9} {\rm\cite{quasi}}
If $M$ is a quasi hemi-slant submanifold of a metallic Riemannian manifold $(\bar M,g,\varphi)$ with the quasi hemi-slant angle $\theta,$ then for all $X, Y \in \Gamma(TM)$ we have
\begin{align*}&g(TP^\theta X,TP^\theta Y) =\cos^2{\theta}[pg(\varphi P^\theta X,P^\theta Y)+qg(P^\theta X,P^\theta Y)],
\\& g(NX,NY)= -\sin^2{\theta}[pg(\varphi P^\theta X,P^\theta Y)+qg(P^\theta X,P^\theta Y)]
\\&\hskip.95in  - [pg(\varphi P^\perp X, P^\perp Y) + qg(P^\perp X, P^\perp Y)],
\end{align*} where $P^\theta$ and $P^\perp$ are the projections of $X \in \Gamma(TM)$ on the distributions $D^\theta$ and $D^\perp$, respectively.
\end{theorem}

\begin{theorem} \label{T:3.10} {\rm\cite{quasi}}
If $M$ is a quasi hemi-slant submanifold of a metallic Riemannian manifold $(\bar M,g,\varphi)$ with the quasi hemi-slant angle $\theta.$ Then $$T^2 P^\theta= \cos^2 {\theta}\,[p\varphi P^\theta +qP^\theta].$$
\end{theorem}

\begin{corollary} \label{C:3.11} {\rm\cite{quasi}} If $M$ is a quasi hemi-slant submanifold of a metallic Riemannian manifold $(\bar M,g,\varphi)$ with the quasi hemi-slant angle $\theta.$ Then $$T^2 P^\theta= \cos^2 {\theta}\, [\varphi +I] P^\theta,$$ where $I$ is the identity mapping on $\Gamma(D^\theta).$
\end{corollary}

\begin{remark}
{\rm Karmakar and Bhattacharyya have established the integrability conditions and certain properties for distributions associated with quasi hemi-slant submanifolds, including an example within a metallic Riemannian manifold, as detailed in \cite{quasi}.}
\end{remark}

\subsection{Pointwise slant submanifolds of metallic Riemannian manifolds}\label{S3.7}

The pointwise slant submanifolds in metallic Riemannian manifolds are studied in \cite{w3} by Hre\c{t}canu and Blaga as follows.

\begin{definition} 
A submanifold $M$ of a metallic Riemannian manifold $(\bar M, g, \varphi)$ is referred to as {\it pointwise slant} if the angle $\theta_x(X)$ between $\varphi X$ and $T_xM$ (known as the Wirtinger angle) remains the same irrespective of the chosen tangent vector $X \in T_xM \setminus \{0\},$ although it may varies with $x \in M$. The Wirtinger angle is a real-valued function $\theta$ (referred to as the Wirtinger function), given by
$$\cos{\theta_x} =\frac{\|TX\|}{\|\varphi X\|},$$
for any $x \in M$ and any $X \in T_xM \setminus \{0\}.$

A pointwise slant submanifold of a metallic Riemannian manifold is called a {\it slant submanifold} if its Wirtinger function $\theta$ is globally constant.
\end{definition}

\begin{proposition} \label{P:3.13} {\rm \cite{w3}}
If $M$ is an isometrically immersed submanifold in the metallic Riemannian manifold $(\bar M, g, \varphi)$, then $M$ is a pointwise slant submanifold if and only if we have $T^2=(\cos^2 {\theta})(pT+qI)$
for some real-valued function $\theta$.
\end{proposition}

\begin{proposition}  \label{P:3.14} {\rm \cite{w3}}
Let $M$ be an isometrically immersed submanifold in the metallic Riemannian manifold $(\bar M, g, \varphi)$. If $M$ is a pointwise slant submanifold with the Wirtinger angle $\theta$, then  we have:
\begin{align*} &g(NX, NY ) = (\sin^2 {\theta})[pg(T X, Y ) + qg(X, Y )]\\& 
tNX = (\sin^2 {\theta})(pT X + qX),\end{align*}
for any $X,Y \in \Gamma(TM)$, where $tV := (\varphi V)^T$ is the tangential component of $\varphi V$ with $V \in \Gamma(TM^\perp)$. \end{proposition}

\begin{remark}
{\rm In \cite{w3}, the study of pointwise bi-slant submanifolds within metallic Riemannian manifolds is conducted, yielding interesting findings.}
\end{remark}

\begin{remark}
{\rm The study of the geometry of submanifolds of co-dimension 2 of locally metallic Riemannian manifolds is done in \cite{ah} and the study of some characterizations of any submanifold of a locally decomposable metallic Riemannian manifold in the case that the co-dimension of the submanifold is greater than or equal to the rank of the set of tangent vector fields of the induced structure on it by the metallic Riemannian structure of the ambient manifold is carried out in \cite{gok} by G\"ok.}
\end{remark}

\section{Submanifolds Immersed in Almost-Complex Metallic \\Manifolds}\label{S4}

In \cite{acm2}, authors have explored the geometric properties of submanifolds within an almost-complex metallic manifold and identified various types of submanifolds, including invariant, anti-invariant, and slant submanifolds.

\begin{proposition} \label{P:4.1} {\rm \cite{acm2}}
Let $M$ be an $n$-dim submanifold of an $(n + s)$-dim ACMSR manifold.
Then, for any $X,Y, Z \in \Gamma(TM)$ and second fundamental form $h$, we have
$$g((\nabla_X h)Y, Z) = g(Y,(\nabla_X h)Z).$$
\end{proposition}

\begin{theorem}  \label{T:4.1} {\rm \cite{acm2}}
Let $M$ be an $n$-dim submanifold of an $(n + s)$-dim ACMSR manifold. If $h$ is parallel, concerning the Levi-Civita connection on $M$, and $V_i$ ($V_i's$) are vector fields on $M$) $(1 \leq i \leq s)$ are linearly independent, then $M$ is totally geodesic.
\end{theorem}

\begin{theorem}  \label{T:4.2} {\rm \cite{acm2}}
Let $M$ be an $n$-dim submanifold of an $(n + s)$-dim ACMSR manifold. $M$ is slant if and only if we have
$$h^2 =\cos^2{\theta}\left(ah-\frac{3}{2}bI\right).$$
\end{theorem}

\section{Warped Product Submanifolds of Metallic Riemannian  Manifolds}\label{S5}

Let $(\bar{M}_1, g_1)$ and $(\bar{M}_2, g_2)$ be two Riemannian manifolds of dimensions $n_1 > 0$ and $n_2 > 0$, respectively. We denote by $P_1$ and $P_2$ the projection maps from the product manifold $\bar{M}_1 \times \bar{M}_2$ onto $\bar{M}_1$ and $\bar{M}_2$, respectively. Let $\bar{\phi} := \phi \circ P_1$ be the lift to $\bar{M}_1 \times \bar{M}_2$ of a smooth function $\phi$ on $\bar{M}_1$. Here, $\bar{M}_1$ is called the base and $\bar{M}_2$ is the fiber of $\bar{M}_1 \times \bar{M}_2$. The unique element $\bar{X}$ of $\Gamma(T(\bar{M}_1 \times \bar{M}_2))$ that is $P_1$-related to $X \in \Gamma(T\bar{M}_1)$ and to the zero vector field on $\bar{M}_2$ will be called the horizontal lift of $X$. Similarly, the unique element $\bar{V}$ of $\Gamma(T(\bar{M}_1 \times \bar{M}_2))$ that is $P_2$-related to $V \in \Gamma(T\bar{M}_2)$ and to the zero vector field on $\bar{M}_1$ will be called the vertical lift of $V$. We denote by $\mathcal{L}(\bar{M}_1)$ the set of all horizontal lifts of vector fields on $\bar{M}_1$ and by $\mathcal{L}(\bar{M}_2)$ the set of all vertical lifts of vector fields on $\bar{M}_2$.

For $f : \bar M_1 \to (0, \infty)$ a smooth function on $M_1$, we consider the Riemannian
metric $g$ on $\bar M :=\bar M_1 \times \bar M_2:$
\begin{equation} \label{g}
g:=P_1^* g_1 +(f \circ P_1)^2 P_2^* g_2.
\end{equation}

\begin{definition} \cite{w2,book17}
 The product manifold of $\bar M_1$ and $\bar M_2$ together with the Riemannian metric $g$ defined by \eqref{g} is called the warped product of $\bar M_1$ and $\bar M_2$ by the warping function $f$.
\end{definition}

Hre\c{t}canu and Blaga investigate the presence of proper warped product submanifolds within metallic Riemannian manifolds in \cite{w1}.

\begin{definition} \cite{w1}
If $M:= M_1 \times_f M_2$ is a warped product submanifold in a metallic Riemannian manifold $(\bar M,g,\varphi )$ such that one of the components $M_i (i \in {1, 2})$ is an invariant submanifold (respectively, anti-invariant submanifold) in $\bar M$ and the other one is a slant submanifold in $\bar M$, with the slant angle $\theta \in [0,\frac{\pi}{2}],$ then
one can call the submanifold $M$ warped product semi-slant (respectively, hemi-slant)
submanifold in the metallic Riemannian manifold $(\bar M,g,\varphi )$.
\end{definition}

\begin{theorem} \label{T:5.1} {\rm \cite{w1}}
Let $M:= M_T \times_f M_\perp$ be a warped product semi-invariant submanifold in a locally metallic Riemannian manifold $(\bar M,g,\varphi )$ (i.e. $M_T$ is invariant and $M_\perp$ is an anti-invariant submanifold in $\bar M$). Then $M:= M_T \times_f M_\perp$ is a non-proper warped product submanifold in $\bar M$ (i.e., the warping function $f$ is constant on the connected components of $M_T$).
\end{theorem}

\begin{theorem}  \label{T:5.2} {\rm \cite{w1}}
Let $M:= M_T \times_f M_\theta$ be a warped product semi-slant submanifold in a locally metallic  Riemannian manifold $(\bar M, g, \varphi)$ (i.e., $M_T$ is invariant and $M_\theta$ is a proper slant submanifold in $\bar M$, with the slant angle $\theta \in (0,\frac{\pi}{2}).$ Then $M := M_T \times_f M_\theta$ is a non proper warped product submanifold in $\bar M$ (i.e. the warping function $f$ is constant on the connected components of $M_T$).
\end{theorem}

\begin{remark}
{\rm The discussion on semi-invariant, semi-slant, and hemi-slant warped product submanifolds within metallic Riemannian manifolds can be found in \cite{w1}, where the authors present various examples of these submanifolds in Euclidean spaces.}
\end{remark}

Hre\c{t}canu and Blaga, in their work \cite{w3}, have explored warped product pointwise bi-slant submanifolds as well as warped product pointwise semi-slant or hemi-slant submanifolds within metallic Riemannian manifolds, deriving several significant results:

\begin{theorem}  \label{T:5.3} {\rm \cite{w3}}
If $M:= M_T \times_f M_\theta$ is a warped product pointwise semi-slant submanifold in a locally metallic Riemannian manifold $(\bar M,g,\varphi)$ with the pointwise slant angle $\theta_x \in (0,\frac{\pi}{2}),$ for $x \in M_\theta,$ then the warping function $f$ is constant on the connected components of $M_T$.
\end{theorem}

\begin{theorem}\label{T:5.4} {\rm \cite{w3}}
If $M:= M_\perp \times_f M_\theta$ (or $M:= M_\theta \times_f M_\perp$) is a warped product pointwise hemi-slant submanifold in a locally metallic Riemannian manifold $(\bar M,g,\varphi)$ with the pointwise slant angle $\theta_x \in (0,\frac{\pi}{2}),$ for $x \in M_\theta,$ then the warping function $f$ is constant on the connected components of $M_\perp$ if and only if we have 
$$A_{NZ} X = A_{NX} Z,$$ where $A$ is the shape operator.
For any $X \in \Gamma(TM_\perp)$ and $Z \in \Gamma(TM_\theta)$ (or $X \in \Gamma(TM_\theta)$ and $Z \in \Gamma(TM_\perp),$ respectively).
\end{theorem}

Bi-warped product submanifolds in some structures of metallic Riemannian manifold are studied in \cite{w5} by Bhunia et al.

\begin{definition} \cite{book17,w4}
Let $\bar M_0, \bar M_1$ and $\bar M_2$ be three Riemannian manifolds and $M = \bar M_0 \times \bar M_1 \times M_2$ be their cartesian product. $P_i : \bar M \to \bar M_i$ is the canonical projection of $\bar M$ onto $\bar M_i$ , where $i \in \{0, 1, 2\}.$ Let $P_i^* : T \bar M \to T\bar M_i$ is the tangent map of $P_i : \bar M \to \bar M_i ,$ where $\Gamma(T\bar M)$ is the Lie algebra of the vector fields of $\bar M$.

If $f_1$ and $f_2$ are two positive real-valued functions on $\bar M_0$, then
\begin{align*}
g(X, Y) &= g(P_0*X, P_0*Y) + (f_1 \circ P_1) 2g(P_1*X, P_1*Y) \\
        &\quad + (f_2 \circ P_2) 2g(P_2*X, P_2*Y).
\end{align*}
$X, Y \in \Gamma(T\bar M)$ defines a Riemannian metric on $\bar M$. This is called the bi-warped product metric.

The product manifold $M = \bar M_0 \times \bar M_1 \times M_2$ furnished by the metric $g$ is called a bi-warped
product manifold and it is denoted by $\bar M_0 \times_{f_1} \bar M_1 \times_{f_2}  \bar M_2$. The functions $f_1$ and $f_2$ are called warping functions. 
\end{definition}

\begin{remark}
    {\rm Examples and some results of bi-warped product submanifolds in some structures of metallic Riemannian manifold are given in \cite{w5}}.
\end{remark}

\section{Lightlike Submanifolds of Metallic Riemannian \\Manifolds}\label{S6}

The main distinction between the theory of lightlike submanifolds and semi-Riemannian submanifolds stems from the fact that, in the former, a portion of the normal vector bundle $TM^\perp$ is contained within the tangent bundle $TM$ of the submanifold $M$ of a semi-Riemannian manifold $\bar{M}$, whereas in the latter, $TM \cap TM^\perp = \{0\}$. Consequently, the fundamental issue in lightlike submanifolds is to substitute the intersecting part with a vector subbundle whose sections are never tangent to $M$. To create a non-intersecting lightlike transversal vector bundle of the tangent bundle, Duggal and Bejancu employed an extrinsic method, while Kupeli utilized an intrinsic method \cite{kupeli}. Since then, numerous researchers have explored the geometry of lightlike hypersurfaces and lightlike submanifolds.

It is a well-established fact that the non-degenerate metric $g$ of a $(m + n)$-dimensional semi-Riemann manifold $\bar M$  is not necessarily induced as a non-degenerate metric on an $m$-dimensional submanifold $M$ of $\bar M$. When the induced metric $g$ is degenerate on $M$ and the rank of $(Rad(TM))$ (where $Rad(TM)$ is subbundle of of null vectors within the tangent bundle 
$TM$) is $r,$ where $1 \leq r \leq m,$ the pair $(M, g)$ is referred as a lightlike submanifold \cite{li}.

\subsection{Lightlike hypersurfaces of metallic semi-Riemannian manifolds}\label{S6.1}

The lightlike hypersurfaces of metallic semi-Riemannian manifolds are introduced by Acet in \cite{acet} and the following results have been introduced:

Let $M$ be a lightlike hypersurface of a metallic semi-Riemannian manifold $(\bar M,g,\varphi)$. For every $X \in \Gamma (TM)$ and $N \in \Gamma ({\rm ltr}(TM)),$ (where $({\rm ltr}(TM))$ is the lightlike transversal bundle of $TM$), $$\varphi X=\varphi' X+\eta (X)N,\;\;\; \varphi N=\xi +\lambda (E)N,$$ where $\varphi' X, \xi \in \Gamma(TM)$ and $\eta,\lambda$ are the 1-forms given by $$\eta (X)=g(X,\varphi E),\;\;\; \lambda (X)=g(X,\varphi N),$$ where $\varphi'$ is a $(1,1)$-tensor field on $M$ and $E$ is a non-zero section of $Rad(TM)$ \cite{acet}.

\begin{lemma} \label{L:6.1} {\rm \cite{acet}}
Let $M$ be a lightlike hypersurface of a locally metallic semi-Riemannian manifold. Then the structure $\varphi'$ is a metallic structure on $M$.
\end{lemma}

\subsubsection{Invariant hypersurface of metallic semi-Riemannian manifolds}\label{S6.1.1}

\begin{definition} \cite{acet}
Let $M$ be a lightlike hypersurface of a locally metallic semi-Riemannian manifold. In that case $M$ is called invariant hypersurface of $\bar M$ if $$\varphi({Rad}(TM))=({Rad}(TM)),\;\;\; \varphi(\operatorname{ltr}(TM))=(\operatorname{ltr}(TM)).$$ 
\end{definition}

\begin{theorem} \label{T:6.1} {\rm \cite{acet}}
Let $M$ be an invariant lightlike hypersurface of a locally metallic semi-Riemannian manifold. Then

{\rm (i)} $h(X, \varphi Y ) = h(\varphi X, Y) = \varphi h (X , Y),$

{\rm (ii)} $h(\varphi X,\varphi Y) = ph(X, \varphi Y) + qh(X, Y).$
\end{theorem}

\subsubsection{Screen semi-invariant hypersurfaces of metallic semi-Riemannian manifolds}\label{S6.1.2}

\begin{definition} \cite{acet}
Let $M$ be a lightlike hypersurface of a locally metallic semi-Riemannian manifold $(\bar M,g,\varphi)$. If
$$\varphi (Rad(TM)) \subset S(TM)\;\;\;{\rm and}\;\;\; 
\varphi ({\rm ltr}(TM)) \subset S(TM)$$ hold, where $S(TM)$ is the screen distribution of the tangent bundle, then $M$ is called a {\it screen semi-invariant hypersurface} of $\bar M$.
\end{definition}

\begin{theorem} \label{T:6.2} {\rm \cite{acet}}
Assume that $M$ is a screen semi-invariant lightlike hypersurface of a locally metallic semi-Riemannian manifold. Then lightlike vector field $\psi$ is parallel on $M$ if and only if we have

{\rm (i)} $M$ is totally geodesic on $\bar M,$

{\rm (ii)} $\gamma = 0.$ ($\gamma$ represents the transversal vector field associated with the screen distribution of the hypersurface $M$.)
\end{theorem}

\begin{theorem} \label{T:6.3} {\rm \cite{acet}}
Let $M$ be a screen semi-invariant lightlike hypersurface of a locally metallic semi-Riemannian manifold. Then lightlike vector field $\Omega$ is parallel on $M$ iff $M$ and $S(TM)$ is totally geodesic on $\bar M$.
\end{theorem}

\begin{remark}
{\rm Results regarding mixed geodesic lightlike hypersurface in metallic semi-Riemannian manifold are discussed in \cite{acet}. Some results on screen semi-invariant lightlike hypersurface of a metallic semi-Riemannian manifold are obtained in \cite{li} by Perktaş et al.}
\end{remark}

\begin{remark}
 {\rm The proof of geodesic and minimal conditions for hypersurfaces of metallic Riemannian manifolds is discussed in \cite{ch}. Choudhary et al. in \cite{ch} investigated $k$-almost Newton-Ricci solitons ($k$-ANRS) embedded in a metallic Riemannian manifold $\bar M$ having the potential function $\psi$ and have explained some applications of metallic Riemannian manifold admitting $k$-almost Newton-Ricci solitons.} 
\end{remark}

\subsection{Invariant lightlike submanifolds of metallic semi-Riemannian \\manifolds}\label{S6.2}

The invariant lightlike submanifold in metallic Riemannian manifold is studied in \cite{li} by Perktaş et al. 

\begin{definition} \cite{li}
Let $(\bar M,g,\varphi)$ be a metallic semi-Riemannian manifold and $(M,g)$ be a lightlike submanifold of $\bar M$. Then $M$ is called an invariant lightlike submanifold of $\bar M$ if the following two conditions are satisfied:
$$\varphi (S(TM)) = S(TM)\;\; {\rm and}\;\;
\varphi (Rad(TM)) = Rad(TM).$$     
\end{definition}

\begin{corollary} \label{C:6.4} {\rm \cite{li}}
Let $(\bar M,g,\varphi)$ be a metallic semi-Riemannian manifold and $(M,g)$ be an invariant lightlike submanifold of $\bar M$. Then the lightlike transversal distribution $ltr(TM)$ is invariant under $\varphi$.
\end{corollary}

\begin{theorem} \label{T:6.5} {\rm \cite{li}}
Let $M$ be an invariant lightlike submanifold of a metallic semi-Riemannian manifold $\bar M$. Then,
the radical distribution $Rad(TM)$ is integrable if and only if we have either $$A^*_{\varphi X}Y=A^*_{\varphi Y}X  \quad  A^*_{X}Y=A^*_{Y}X$$
or
$$A^*_{\varphi X}Y-A^*_{\varphi Y}X=p(A^*_{X}Y=A^*_{Y}X)$$
for any $X,Y \in \Gamma(Rad(TM))$ and any $Z \in \Gamma(S(TM)),$ where $A^*$ denotes the shape operator of the distributions $(S(TN))$ and $Rad(TN)$.
\end{theorem}

\begin{theorem}\label{T:6.6} {\rm \cite{li}}
Let $M$ be an invariant lightlike submanifold of a metallic semi-Riemannian manifold $\bar M$. Then,
the induced connection $\nabla$ on $M$ is a metric connection if and only if we have $$A^*_{\varphi \xi}U=pA^*_{\xi}U$$
for  any $U \in \Gamma(TM)$ and any $\xi \in \Gamma(Rad(TM)).$
\end{theorem}

\begin{remark}
    {\rm Results concerning totally geodesic foliation of an invariant lightlike submanifold of a metallic semi-Riemannian manifold are discussed in \cite{li}.}
\end{remark}

The study of the geometry of the semi-invariant lightlike submanifold of a metallic semi-Riemannian manifold is conducted in \cite{s1} by Kaur et al. and some conditions for the integrability of distributions are discussed.

\begin{definition} \cite{s1}
Let $M$ be a lightlike submanifold of a metallic semi-Riemannian manifold $(\bar M,g,\varphi)$. Then, $M$
is a semi-invariant lightlike submanifold, if the following three conditions are satisfied:
$$\varphi (Rad(TM)) \subseteq S(TM),\;\; \varphi (ltr(TM)) \subseteq S(TM)\;\; {\rm and}\;\; \varphi S(TM^\perp) \subseteq S(TM).$$
\end{definition}

\begin{theorem} \label{T:6.7} {\rm \cite{s1}}
Let $M$ be a semi-invariant lightlike submanifold of a metallic semi-Riemannian manifold $(\bar M,g,\varphi)$. Then $D$ is integrable if and only if $$h^l (\varphi X,\varphi Y)=ph^l (\varphi X, Y)+qh^l (X, Y),$$ $$h^s (\varphi X,\varphi Y)=ph^s (\varphi X, Y)+qh^s (X, Y)$$ hold for $X, Y \in \Gamma(D)$, where $h^l$ and $h^s$ denote the second fundamental forms associated with the lightlike and screen distributions of submanifold $M$, respectively. 
\end{theorem}

\begin{theorem} \label{T:6.8} {\rm \cite{s1}}
Let $M$ be a semi-invariant lightlike submanifold of a metallic semi-Riemannian manifold $(\bar M,g,\varphi)$. If $D^\theta$ is integrable, then  leaves of $D^\theta$ have a metallic structure.
\end{theorem}

\begin{remark}
{\rm Moreover, the authors in \cite{s1} have investigated totally geodesic and mixed geodesic distributions of semi-invariant lightlike submanifolds.}
\end{remark}

The structure of invariant and screen
semi-invariant lightlike submanifolds of a metallic semi-Riemannian manifold with a quarter symmetric non-metric connection is introduced in \cite{s2}.

\begin{theorem} \label{T:6.9} {\rm \cite{s2}}
For an invariant lightlike submanifold $M$ of a metallic semi-Rieman\-nian manifold $\bar M$ with a quarter symmetric non-metric connection $\mathcal{\bar D}$, the radical distribution is
integrable if and only if we have $$A^*_{\varphi U}V-pA^*_{U}V=A^*_{\varphi V}U-pA^*_{V}U,$$ for any $U,V \in \Gamma(Rad(TM))$ and $Z \in \Gamma(S(TM)).$
\end{theorem}

\begin{theorem} \label{T:6.10} {\rm \cite{s2}}
Let $M$ be an invariant lightlike submanifold of a metallic semi-Riemannian manifold $\bar M$ with a quarter symmetric non-metric connection $\mathcal{\bar D}$. Then the screen 
distribution is integrable if and only if 
$$A^* (V,\varphi U)+pA^* (U,V)=A^* (U,\varphi V)+pA^* (V,U)$$ holds for any $U,V \in \Gamma(S(TM))$ and any $N \in \Gamma({\rm ltr}(TM)).$
\end{theorem}

The results on the structure of screen semi-invariant lightlike submanifolds of a metallic 
semi-Riemannian manifold with a quarter symmetric non-metric connection are the following.

\begin{proposition} \label{P:6.1} {\rm \cite{s2}}
Let $M$ be a screen semi-invariant lightlike submanifold of a metallic semi-Riemannian manifold. Then $M$ is an invariant lightlike
submanifold of $\bar M$ if and only if $D^\perp = \{0\},$ where $D^\perp$ is orthogonal complementary distribution of $D^\theta .$
\end{proposition}

\begin{proposition}\label{P:6.2} {\rm \cite{s2}}
If a screen semi-invariant lightlike submanifold of a metallic semi-Riemannian manifold $\bar M$ is isotropic or totally lightlike, then it is an invariant lightlike submanifold of $\bar M.$
\end{proposition}

\begin{theorem} \label{T:6.11} {\rm \cite{s2}}
Let $M$ be a screen semi-invariant lightlike submanifold of metallic semi-Riemannian manifold $\bar M$ with a quarter symmetric non-metric connection $\mathcal{\bar D}$. Then the necessary and sufficient condition for $D^\theta$ to be integrable is that
\begin{align*}&h^s (U,\varphi V)=h^s(V,\varphi U),  
\\&h^*(U,\varphi V)+ph^*(V,U)=h^*(V,\varphi U)+ph^* (U,V) \end{align*}
for any $U,V \in \Gamma(D^\theta), \,Z \in \Gamma(D^\perp)$ and $N \in \Gamma({\rm ltr}(TM))$.
\end{theorem}

\begin{remark}
{\rm \cite{s2} contains the results regarding totally geodesic foliations for invariant and screen semi-invariant submanifolds.}
\end{remark}

\subsection{Transversal lightlike submanifolds of metallic semi-Riemannian \\manifolds}\label{S6.3}

The study of the geometry of transversal lightlike submanifolds and radical transversal lightlike submanifolds of metallic semi-Riemannian manifolds is proposed in \cite{tr}. The authors in \cite{tr} explored the geometry of distributions and established the essential and adequate conditions for the induced connections in these manifolds to qualify as metric connections. Furthermore, they provided a characterization of transversal lightlike submanifolds of metallic semi-Riemannian manifolds.

Results on radical transversal lightlike submanifolds of a metallic semi-Riemannian manifold are as follows:

\begin{definition} \cite{tr}
Let $(M, g, S(TM), S(T M^\perp))$ be a lightlike submanifold of a metallic semi-Riemannian manifold $(\bar M,g,\varphi).$ Then  $M$ is called a {\it radical transversal lightlike submanifold} if the following two conditions are satisfied:
$$\varphi Rad (TM)= {\rm ltr}(TM)\;\; {\rm and}\;\;\varphi (S(TM)) = S(TM).$$
\end{definition}

\begin{proposition} \label{P:6.3} {\rm \cite{tr}}
Let $\bar M$ be a metallic semi-Riemannian manifold. In this case, there is no 1-radical transversal lightlike submanifold of $\bar M$.     
\end{proposition}

\begin{theorem}  \label{T:6.12} {\rm \cite{tr}} Let $M$ be a radical transversal lightlike submanifold of a metallic semi-Riemannian manifold $\bar M$. In this case,  $S(TM^\perp)$ is invariant under $\varphi$. 
\end{theorem}

\begin{theorem}  \label{T:6.13} {\rm \cite{tr}}
Let $M$ be a radical transversal lightlike submanifold of a locally metallic semi-Riemannian manifold. Then the induced connection $\nabla$ on $M$ is a metric connection if and only if there is no component of $A_{\varphi \xi}W$ in $\Gamma (S(TM))$
for any $W \in \Gamma(TM)$ and $\xi \in \Gamma(Rad (TM)).$
\end{theorem}

\begin{theorem}  \label{T:6.14} {\rm \cite{tr}}
Let $M$ be a radical transversal lightlike submanifold of a locally metallic semi-Riemannian manifold $\bar M$. In this case, the screen distribution is integrable if and only if we have
$h^l (U, SW) = h^l (W, SU)$
for any $W, U \in \Gamma(S(TM)).$
\end{theorem}

\begin{definition} \cite{tr}
Let $(M, g, S(TM), S(TM^\perp))$ be a lightlike submanifold of a metallic semi-Riemannian manifold $(\bar M,g,\varphi)$. Then $M$ is called a {\it transversal lightlike submanifold} if the following two conditions are satisfied:
$$\varphi Rad(TM)=\operatorname{ltr}(TM)\;\; {\rm and}\;\; \varphi(S(TM))\subseteq S(TM^\perp).$$
\end{definition}

The following results on transversal lightlike submanifolds of metallic semi-Riemannian manifolds are obtained in \cite{tr}:

\begin{proposition}  \label{P:6.4} {\rm \cite{tr}}
Let $M$ be a transversal lightlike submanifold of a locally metallic semi-Riemannian manifold $\bar M$. In
this case, the distribution $\mu$ is invariant according to $\varphi$, where $\mu$ denotes the orthogonal complement subbundle to $\varphi (S(TM))$ in $S(TM^\perp)$.
\end{proposition}

\begin{proposition}  \label{P:6.5} {\rm \cite{tr}}
There does not exist a $1$-lightlike transversal lightlike submanifold of a locally metallic semi-Riemannian manifold.    
\end{proposition}

\begin{corollary}  \label{C:6.15} {\rm \cite{tr}}
Let $M$ be a transversal lightlike submanifold of a locally metallic semi-Riemannian manifold $\bar M$. Then:

{\rm (i)} $\dim (Rad (TM)) \geq 2.$

{\rm (ii)} The transversal lightlike submanifold of $3$-dimensional is $2$-lightlike.
\end{corollary}

\begin{theorem}  \label{T:6.16} {\rm \cite{tr}}
Let $M$ be a transversal lightlike submanifold of a locally metallic semi-Riemannian manifold $\bar M$. Then
the radical distribution is integrable if and only if
$$D^s (U,LW)=D^s (W,LU),$$
for $U, W \in \Gamma(Rad (TM))$, where $D^s$ denotes the screen distribution and $L$ represents the structure endomorphism associated with the metallic structure.
\end{theorem}

\begin{remark}
{\rm The conditions of a totally geodesic foliation of such submanifolds and necessary and sufficient conditions for the induced connection on these manifolds to be a metric
connection in metallic semi-Riemannian manifold is obtained in \cite{tr}.}
\end{remark}

In \cite{tr2}, Erdoğan et al. presented the concept of screen transversal lightlike submanifolds within metallic semi-Riemannian manifolds, along with its subclasses: screen transversal anti-invariant, radical screen transversal, and isotropic screen transversal lightlike submanifolds.

\begin{definition} \cite{tr2}
Let $M$ be a lightlike submanifold of a metallic semi-Riemannian manifold $(\bar M,g,\varphi).$ If
$\varphi Rad (TM) \subset S(TM^\perp)$ holds,
then $M$ is called a {\it screen transversal lightlike submanifold} of a metallic semi-Riemannian manifold.
\end{definition}

\begin{definition} \cite{tr2}
Let $M$ be a screen transversal lightlike submanifold of a metallic semi-Riemannian manifold $(\bar M,g,\varphi).$

{\rm (i)} If $\varphi S(TM) \subset S(TM^\perp),$ then $M$ is called a {\it screen transversal anti-invariant lightlike submanifold} of $(\bar M,g,\varphi).$

{\rm (ii)} If $\varphi S(TM) = S(TM),$ then we say that $M$ is a {\it radical screen transversal lightlike submanifold}of $(\bar M,g,\varphi).$
\end{definition}

\begin{proposition}  \label{P:6.6} {\rm \cite{tr2}}
Let $M$ be a screen transversal anti-invariant lightlike submanifold of a metallic semi-Riemannian manifold $(\bar M,g,\varphi).$ Then the distribution D is invariant for $\varphi$.
\end{proposition}

\begin{proposition}  \label{P:6.7} {\rm \cite{tr2}}
Let $M$ be a screen transversal anti-invariant lightlike submanifold of a metallic semi-Riemannian manifold $(\bar M,g,\varphi).$ Then there do not exist co-isotropic and totally screen transversal types of such lightlike submanifolds.
\end{proposition}

\begin{proposition}  \label{P:6.8} {\rm \cite{tr2}}
Let $M$ be a radical screen transversal lightlike submanifold of a locally metallic semi-Riemannian manifold $(\bar M,g,\varphi).$ Then the distribution D is invariant concerning $\varphi$.
\end{proposition}

\begin{remark}
{\rm The structure of distributions relevant to defining these submanifolds and the criteria for the induced connection to qualify as a metric connection. Additionally, the necessary and sufficient condition for an isotropic screen transversal lightlike submanifold to be totally geodesic is examined in \cite{tr2}.}
\end{remark}

\begin{remark}
{\rm In \cite{um}, Shankar and Yadav have examined the geometry of totally umbilical screen-transversal lightlike submanifolds. They investigated two categories: totally umbilical radical screen-transversal lightlike submanifolds and totally umbilical screen-transversal anti-invariant lightlike submanifolds. The authors derived the necessary and sufficient conditions for the integrability of the distributions and for the induced connection to be a Levi-Civita or metric connection on these lightlike submanifolds.}
\end{remark}

\subsection{Half-lightlike submanifolds of metallic semi-Riemannian manifolds}\label{S6.4}

A lightlike submanifold of co-dimension 2 of a semi-Riemannian manifold is called a half-lightlike
submanifold if the mapping defining the radical distribution has rank 1 \cite{h1,h2}. Screen conformal half-lightlike submanifolds of semi-Riemannian manifolds are presented in \cite{h2}.

The study of half-lightlike submanifolds of a semi-Riemannian manifold endowed with a metallic structure is done by Acet et al. in \cite{half}.

\begin{definition} \cite{half}
Assume that $M$ is called a {\it half-lightlike submanifold} of a metallic semi-Riemannian manifold $(\bar M,g,\varphi).$ If we have
$$\varphi (Rad(TM)) \subset S(TM),\;\; \varphi ({\rm ltr}(TM)) \subset S(TM),\;\;{\rm and}\;\;  \varphi (S(TM^\perp)) \subset S(TM),$$
then $M$ is called a {\it screen semi-invariant half-lightlike submanifold} of $\bar M$. 
\end{definition}

\begin{theorem} \label{T:6.17} {\rm \cite{hemi}}
Assume that $M$ is a half-lightlike submanifold of a metallic semi-Riemannian manifold $(\bar M,g,\varphi).$
Then $\phi$ is a metallic structure on $D^\theta$.
\end{theorem}

\begin{theorem} \label{T:6.18} {\rm\cite{half}}
Assume that $M$ is a half-lightlike submanifold of a metallic semi-Riemannian manifold $(\bar M,g,\varphi).$  If the distribution $D^\perp$ is parallel then $D^\perp$ is totally geodesic on $M$.
\end{theorem}

\subsection{Slant lightlike submanifolds of metallic semi-Riemannian manifolds}\label{S6.6}

Lone and Harry have introduced in \cite{lone}  the notion of screen slant lightlike submanifolds of metallic semi-Riemannian manifolds and obtained the following results:

\begin{definition} \cite{lone}
 Let $M$ be a $2q$-lightlike submanifold (where $q$ is an integer indicating the dimension of the radical bundle of the lightlike submanifold $M$ within the metallic semi-Riemannian manifold) of metallic semi-Riemannian manifold $\bar M$ of index $2q < \dim(M).$ Then we say that $M$ is a {\it screen slant lightlike submanifold} of $\bar M$ if the following two conditions are satisfied: 

{\rm (i)} $Rad(TM)$ is invariant with respect to $\varphi$, i.e., $\varphi (Rad(TM)) = Rad(TM).$

{\rm (ii)} $S(TM)$ is slant with $\theta (\neq 0),$ i.e. for each $x \in M$ and each non-zero vector $X \in \Gamma (S(TM)),$ the angle $\theta$ between $\varphi X_1$ and the vector subspace $S(TM)$ is a non-zero constant, which is independent of the choice of $x \in M$ and $X_1 \in \Gamma (S(TM)).$

Furthermore, the constant angle $\theta$ is called the slant angle of distribution $S(TM)$. A  screen lightlike submanifold is said to be proper if it is neither invariant $(\theta = 0)$ nor screen
real $(\theta = \frac{\pi}{2})$.
\end{definition}

\begin{theorem} \label{T:6.21} {\rm \cite{lone}}
Let $M$ be a $2q$-lightlike submanifold of a metallic semi-Riemannian manifold. Then $M$ is a screen slant lightlike submanifold of $\bar M$ if and only if

{\rm (i)} ${\rm ltr} (TM)$ is invariant with respect to $\varphi$,

{\rm (ii)} there exists a constant $c \in [0, 1)$ such that $P^2 X_1 = c (pPX_1 + q X_1),$ for any $X_1 \in \Gamma(S(TM))$ and $PX_1$ is tangential part of $\varphi X_1.$ Moreover, in this case $c = \cos^2 {\theta}$ and $\theta$ is the slant angle of $S(TM).$
\end{theorem}

\begin{remark}
{\rm The necessary and sufficient conditions for the induced connection to be a metric connection. Moreover, investigation of some equivalent conditions for integrability of such submanifolds is obtained in \cite{lone}.}
\end{remark}

\subsection{Warped product lightlike submanifolds of metallic semi-Riemannian manifolds}\label{S6.7}

The investigation of whether the screen real lightlike submanifolds of metallic semi-Riemannian manifolds are warped product lightlike submanifolds or not is done by Shanker and Yadav  in \cite{wl}.

\begin{theorem} \label{T:6.22} {\rm \cite{wl}}
Let $M=M_1 \times_{f} M_2$ be a warped product lightlike submanifold. Then, for any $\xi \in \Gamma (Rad(TM))$ and $U \in \Gamma (S(TM))$, we have $\nabla_\xi U \in \Gamma (S(TM)).$
\end{theorem}

\begin{theorem}
Let $(M,g, S(TM))$ be an irrotational screen-real $r$-lightlike submanifold of a metallic semi-Riemannian manifold, then the induced connection is a metric connection.
\end{theorem}

\begin{theorem}  \label{T:6.23} {\rm \cite{wl}}
There does not exist any class of irrotational screen-real $r$-lightlike submanifolds that can be written in the form of warped product lightlike submanifolds.
\end{theorem}

\begin{remark}
{\rm Some results on a screen-real lightlike submanifold of a metallic semi-Riemannian manifold are also obtained in \cite{wl}.}
\end{remark}

\section{Inequalities in Metallic Riemannian Manifolds}\label{S7}

\subsection{Inequalities involving $\delta$-Casorati curvature in metallic Riemannian manifolds}\label{S7.1}

The Casorati curvature for surfaces in a Euclidean $3$-space $\mathbb E^3$ was introduced by Casorati in 1890 \cite{11}. The Casorati curvature was used in \cite{11} in place of the conventional Gauss curvature. The Casorati curvature $C$ of a submanifold in a Riemannian manifold is defined generally as the normalized squared norm of the second fundamental form. Decu et al. developed normalized Casorati curvatures $\delta_C(n - 1)$ and $\hat{\delta}_C (n - 1)$ in 2007 (see \cite{12}) in keeping with the spirit of $\delta$-invariants. They extended normalized Casorati curvatures in 2008 to generalized normalized $\delta$-Casorati curvatures $\delta_C(r; n-1)$ and $\hat{\delta}_C(r; n-1)$ in \cite{13}. Which laid the groundwork for the definition of optimal inequalities for submanifolds in various ambient spaces using Casorati curvatures.

Let us assume that $(\bar M, g)$ is an $m$-dimensional Riemannian manifold and let $(M, g)$ be an $n$-dimensional Riemannian submanifold isometrically immersed into $\bar M$. Let us denote by $\bar{\nabla}$ the Levi-Civita connection on $\bar M$ and by $\nabla$ the covariant differentiation induced on $M$. Let $h$ be the second fundamental form of $M$ and let $\nabla^{\perp}$ be the connection in the normal bundle. The Gauss and Weingarten formulas are given respectively by
$$\bar{\nabla}_X Y=\nabla_X Y+h(X, Y),\;\; \bar{\nabla}_X \xi=-A_\xi X+\nabla_X^{\perp} \xi $$
for all $X, Y \in \Gamma(TM)$ and $\xi \in \Gamma\left(TM^{\perp}\right)$, where $A_\xi$ denotes the shape operator of $M$ with respect to $\xi$. We also have the following relation between $A_\xi$ and $h$
$$
g\left(A_\xi X, Y\right)=g(h(X, Y), \xi),
$$
for all $X, Y \in \Gamma(TM)$ and $\xi \in \Gamma\left(TM^{\perp}\right)$.
The Gauss equation is \cite{15}
$$
\bar{R}(X, Y, Z, W)=R(X, Y, Z, W)-g(h(X, W), h(Y, Z))+g(h(X, Z), h(Y, W)),
$$
for all $X, Y, Z, W \in \Gamma(TM)$.
Let us denote by $\left\{E_1, \ldots, E_n\right\}$ a local orthonormal tangent frame of the tangent bundle $TM$ of $M$ and by $\left\{E_{n+1}, \ldots, E_m\right\}$ a local orthonormal normal frame of the normal bundle $TM^{\perp}$ of $M$ in $\bar M$. Then the scalar curvature $\tau$ is expressed as
$$
\tau=\sum_{1 \leq i<j \leq n} R\left(E_i, E_j, E_j, E_i\right).
$$
The normalized scalar curvature $\rho$ and the mean curvature $\mathcal{H}$ of $M$ are defined by
$$\rho=\frac{2 \tau}{n(n-1)}\;\; {\rm and}\;\;
\mathcal{H}=\sum_{i=1}^n \frac{1}{n} h\left(E_i, E_i\right) .$$

For our convenience, we put $h_{i j}^r=g\left(h\left(E_i, E_j\right), E_r\right)$, where $i, j=\{1, \ldots, n\}$ and $r=\{n+1, \ldots, m\}$ so that the squared norms of the mean curvature vector $\mathcal{H}$ and of the second fundamental form $h$ can be written as
$$\|\mathcal{H}\|^2=\frac{1}{n^2} \sum_{r=n+1}^m\left(\sum_{i=1}^n h_{i i}^r\right)^2,
\;
\;\; \|h\|^2=\sum_{r=n+1}^m \sum_{i, j=1}^n\left(h_{i j}^r\right)^2 .$$
We recall that the Casorati curvature $C$ of $M$ is given by
$C=\frac{1}{n}\|h\|^2 .$

Let $\left\{E_1, \ldots, E_t\right\}$ be an orthonormal basis of the $t$-dimensional subspace $\mathcal{L}$ of $TM$ with $ t \geq 2$. Then the scalar curvature of $\mathcal{L}$ is
$$
\tau(\mathcal{L})=\sum_{1 \leq i<j \leq t} R\left(E_i, E_j, E_j, E_i\right)
$$
and the Casorati curvature of $\mathcal{L}$ is given by
$$
C(\mathcal{L})=\frac{1}{t} \sum_{r=n+1}^m \sum_{i, j=1}^t\left(h_{i j}^r\right)^2 .
$$

The normalized $\delta$-Casorati curvatures $\delta_C (n-1)$ and $\widehat{\delta}_C (n-1)$ are defined respectively by (see \cite{13})
$$
\left[\delta_C(n-1)\right]_p=\frac{1}{2} C_p+\frac{n+1}{2 n} \inf \left\{C(\mathcal{L}) \mid \mathcal{L}: \text { a hyperplane of } T_p M\right\}
$$
and
$$
\left[\widehat{\delta}_C(n-1)\right]_p=2 C_p-\frac{2 n-1}{2 n} \sup \left\{C(\mathcal{L}) \mid \mathcal{L}: \text { a hyperplane of } T_p M\right\}.
$$
And the generalized normalized $\delta$-Casorati curvatures $\delta_C(r ; n-1)$ and $\widehat{\delta}_C(r ; n-1)$ of the submanifold $M^n$ are defined for any positive real number $r \neq n(n-1)$ as
$$
\begin{aligned}
\left[\delta_C(r ; n-1)\right]_p=r C_p +\frac{(n \! - \!1)(n\! + \!r)\left(n^2\! - \!n- \! r\right)}{r n} \inf \left\{C(\mathcal{L}) \mid \mathcal{L}: \text {a hyperplane of } T_p M\right\}
\end{aligned}
$$
if $0<r<n(n-1)$, and
$$
\begin{aligned}
\left[\widehat{\delta}_C(r ; n-1)\right]_p\! =r C_p +\frac{(n \! - \!1)(n\! + \!r)\left(n^2\! - \!n- \! r\right)}{r n} \sup \left\{C(\mathcal{L}) \mid \mathcal{L}: \text {a hyperplane of } T_p M\right\}
\end{aligned}
$$
for $r>n(n-1)$.

In \cite{14}, Choudhary and Blaga examined sharp inequalities for a slant submanifold using generalized normalized $\delta$-Casorati curvatures in metallic Riemannian space forms $(\bar M = M_p (c_p) \times M_q (c_q), g,\varphi)$ and drawn the following findings:
 
\begin{theorem} \label{T:7.1} {\rm \cite{14}}
Consider an $n$-dimensional slant submanifold $M$ of an $m$-dimensional locally metallic product space form $(\bar{M} = M_p (c_p) \times M_q (c_q), g,\varphi)$. Then

{\rm (i)} The generalized normalized $\delta$-Casorati curvature $\delta_C(r; n-1)$ satisfies
\begin{equation} \label{5}
\begin{aligned}
\rho \leq & \frac{\delta_C(r ; n-1)}{n(n-1)} + \frac{1}{2(p^2 + 4q)} (c_p +c_q ) \left\{p^2 +2q \right. \\
& \left. +\frac{2}{n(n-1)} \left[ \operatorname{tr}^2 \varphi - (p \cdot \operatorname{tr} T+nq) \cos^2{\theta} \right] -\frac{2p}{n} \operatorname{tr}\varphi \right\} \\
& + \frac{1}{2 \sqrt{(p^2 + 4q)}} (c_p -c_q ) \left(\frac{2}{n} \operatorname{tr} \varphi -p\right)
\end{aligned}
\end{equation}
for any real number $0<r<n(n-1);$

{\rm (ii)} The generalized normalized $\delta$-Casorati curvature $\widehat{\delta}_C(r ; n-1)$ satisfies
\begin{equation} \label{6}
\begin{aligned}
\rho \leq & \frac{\widehat{\delta}_C(r ; n-1)}{n(n-1)} + \frac{1}{2(p^2 + 4q)} (c_p +c_q ) \left\{p^2 +2q \right. \\
& \left. +\frac{2}{n(n-1)} \left[ \operatorname{tr}^2 \varphi - (p \cdot \operatorname{tr} T+nq) \cos^2{\theta} \right] -\frac{2p}{n} \operatorname{tr}\varphi \right\} \\
& + \frac{1}{2 \sqrt{(p^2 + 4q)}} (c_p -c_q ) \left(\frac{2}{n} \operatorname{tr} \varphi -p\right)
\end{aligned}
\end{equation}
for any real number $r>n(n-1).$

In addition, the equalities in \eqref{5} and \eqref{6} hold if and only if the submanifold $M$ is invariantly quasi-umbilical with trivial normal connection in $\bar M$ , such that the shape operators $A_r, r \in {n + 1,...,m}$ with respect to some orthonormal tangent frame $\{E_1,...,E_n\}$ and orthonormal normal frame $\{E_n+1,...,E_m\}$
take the following forms:

\begin{equation}
A_{n+1}=\left(\begin{array}{cccccc}
b & 0 & 0 & \ldots & 0 & 0 \\
0 & b & 0 & \ldots & 0 & 0 \\
0 & 0 & b & \ldots & 0 & 0 \\
\vdots & \vdots & \vdots & \ddots & \vdots & \vdots \\
0 & 0 & 0 & \ldots & b & 0 \\
0 & 0 & 0 & \ldots & 0 & \frac{n(n-1)}{r} b
\end{array}\right), \quad A_{n+2}=\cdots=A_m=0.
\end{equation}
\end{theorem}

For normalized $\delta$-Casorati curvature;

\begin{theorem} \label{T:7.2} {\rm \cite{14}}
Consider an $n$-dimensional slant submanifold $M$ of an $m$-dimensional locally metallic product space form $(\bar{M} = M_p (c_p) \times M_q (c_q), g,\varphi)$. Then

{\rm (i)} The normalized $\delta$-Casorati curvature $\delta_C(n-1)$ satisfies
\begin{equation} \label{8}
\begin{aligned}
\rho \leq &\; \delta_C(n-1) + \frac{1}{2(p^2 + 4q)} (c_p +c_q ) \Bigg\{p^2 +2q \\
&  +\frac{2}{n(n-1)} \left[ \operatorname{tr}^2 \varphi - (p \cdot \operatorname{tr} T+nq) \cos^2{\theta} \right] -\frac{2p}{n} \operatorname{tr}\varphi \Bigg\} \\
& + \frac{1}{2 \sqrt{(p^2 + 4q)}} (c_p -c_q ) \left(\frac{2}{n} \operatorname{tr} \varphi -p\right).
\end{aligned}
\end{equation}

{\rm (ii)} The normalized $\delta$-Casorati curvature $\widehat{\delta}_C (n-1)$ satisfies

\begin{equation} \label{9}
\begin{aligned}
\rho \leq &\; \widehat{\delta}_C(n-1)+ \frac{1}{2(p^2 + 4q)} (c_p +c_q ) \Bigg\{p^2 +2q \\
&  +\frac{2}{n(n-1)} \left[ \operatorname{tr}^2 \varphi - (p \cdot \operatorname{tr} T+nq) \cos^2{\theta} \right] -\frac{2p}{n} \operatorname{tr}\varphi \Bigg\} \\
& + \frac{1}{2 \sqrt{(p^2 + 4q)}} (c_p -c_q ) \left(\frac{2}{n} \operatorname{tr} \varphi -p\right),
\end{aligned}
\end{equation}
for any real number $r>n(n-1).$

Moreover, the equalities hold in \eqref{8} and \eqref{9} if and only if $M$ is an invariantly quasi-umbilical submanifold with trivial normal connection in $\bar M$ such that with some orthonormal tangent frame $\{E_1,...,E_n\}$ of the tangent bundle $TM$ of $M$ and orthonormal normal frame $\{E_n+1,...,E_m\}$ of the normal bundle $TM^\perp$ of $M$ in $\bar M$, the shape operators $A_r$ satisfy
\begin{equation}
A_{n+1}=\left(\begin{array}{cccccc}
b & 0 & 0 & \ldots & 0 & 0 \\
0 & b & 0 & \ldots & 0 & 0 \\
0 & 0 & b & \ldots & 0 & 0 \\
\vdots & \vdots & \vdots & \ddots & \vdots & \vdots \\
0 & 0 & 0 & \ldots & b & 0 \\
0 & 0 & 0 & \ldots & 0 & 2 b
\end{array}\right), \quad A_{n+2}=\cdots=A_m=0.
\end{equation}
\end{theorem}

\begin{remark}
{\rm The sharp inequality for invariant and anti-invariant submanifolds in metallic Riemannian space forms was also established by Choudhary and Blaga in \cite{14}, where they further characterized the submanifolds for which the equality holds.}
\end{remark}

\subsection{Chen-type inequality in metallic Riemannian manifolds}\label{S7.2}

One of the fundamental connections between the main extrinsic and intrinsic invariants of a submanifold was examined by the first in 1993 and as a result, he introduced a new Riemannian invariant, $\delta(2)$ \cite{16}. The mathematical description of this invariant is $\delta(2) = \tau - \inf K$, where $K$ signifies the sectional curvature and $\tau$ is the scalar curvature. Chen proved the subsequent inequality using this invariant:
\begin{equation}
    \delta(2) \leq \frac{n-2}{2} \sqrt{n(n-1)} \|\mathcal{H}\|^2 + (n + 1)c,
\end{equation}
for any submanifold $M$ in a real space form with constant sectional curvature $c$ (with $n = \dim M \geq 3$), where  $\|\mathcal{H}\|^2$ denotes the squared mean curvature of $M$. This inequality is referred to as the first Chen inequality.

Inspired by this, many mathematicians considered Chen-type inequalities in metallic Riemannian manifolds. In \cite{17}, Choudhary and Uddin obtained the following results on Chen-type Inequalities for slant submanifolds in metallic Riemannian space forms.

\begin{theorem} \label{T:7.3} {\rm \cite{17}}
Let $M$ be any proper $\theta$-slant submanifold that is isometrically immersed in $(\bar{M} = M_p (c_p) \times M_q (c_q), g,\varphi)$. Then the following inequality follows:
\begin{equation*} \label{12}
\begin{aligned}
    \delta_M \leq & \frac{(n-2)}{2} \left[\frac{4}{(n-1)} \frac{n^2}{4} \|\mathcal{H}\|^2 +\frac{1}{2(p^2 +4q)}(c_p +c_q )\{(p^2 +2q)(n+1)-2p\operatorname{tr}(\varphi)\}\right] \\
    & +\frac{1}{2(p^2 +4q)} (c_p +c_q )  \left[(p\operatorname{tr}(T)-4q)\cos^2{\theta} -\operatorname{tr}^2 (\varphi)\right]\\
    & + \frac{1}{4\sqrt{p^2 +4q}} (c_p -c_q )(n-2) \left[2 \operatorname{tr} (\varphi)-pn+p\right].
\end{aligned}
\end{equation*}
\end{theorem}

\begin{corollary}  \label{C:7.4} {\rm \cite{17}}
For a $\varphi$-invariant submanifold $M^n$ immersed in $\bar M$, the following inequality holds true:
\begin{equation*} 
\begin{aligned}
    \delta_M \leq & \,\frac{(n-2)}{2} \left[\frac{4}{(n-1)} \frac{n^2}{4} \|\mathcal{H}\|^2 +\frac{1}{2(p^2 +4q)}(c_p +c_q )\{(p^2 +2q)(n+1)-2p\operatorname{tr}(\varphi)\}\right] \\
    & +\frac{1}{2(p^2 +4q)} (c_p +c_q )  \left[(p\operatorname{tr}(T)+4q) -\operatorname{tr}^2 (\varphi)\right]\\
    & + \frac{1}{4\sqrt{p^2 +4q}} (c_p -c_q )(n-2) \left[2 \operatorname{tr} (\varphi)-pn-p\right].
\end{aligned}
\end{equation*}
\end{corollary}

\begin{corollary}  \label{C:7.5} {\rm \cite{17}}
A $\varphi$-anti-invariant submanifold $M^n$ immersed in $\bar M$ satisfies the following inequality.
\begin{equation*} 
\begin{aligned}
    \delta_M  \leq &\, \frac{(n-2)}{2} \Bigg[\frac{4}{(n-1)} \frac{n^2}{4} \|\mathcal{H}\|^2 +\frac{1}{2(p^2 +4q)}(c_p +c_q )\{(p^2 +2q)(n+1)-2p\operatorname{tr}(\varphi)\} 
    \\&\hskip.2in -\operatorname{tr}^2 (\varphi) \Bigg] + \frac{1}{4\sqrt{p^2 +4q}} (c_p -c_q )(n-2) \left[2 \operatorname{tr} (\varphi)-pn-p\right].
\end{aligned}
\end{equation*}
\end{corollary}
\begin{remark}
    {\rm In reference \cite{17}, Choudhary and Uddin further demonstrated instances of inequalities for the Ricci curvature tensor in $\theta$-slant, $\varphi$-invariant, and $\varphi$-anti-invariant submanifolds within metallic Riemannian space forms.}
\end{remark}

\subsection{Wintgen inequality in metallic Riemannian manifolds}\label{S7.3}

In 1979, P. Wintgen \cite{18} established a fundamental inequality for the surface $\mathbb M^2$ within the Euclidean 4-space $\mathbb E^4$, commonly known as the Wintgen inequality, which incorporates both intrinsic and extrinsic invariants. He demonstrated that the intrinsic Gaussian curvature $K$ and the extrinsic normal curvature $K^\perp$ of $\mathbb M^2$ in $\mathbb E^4$ fulfill $$K+|K^\perp| \leq \|\mathcal{H}\|^2,$$ where $\|\mathcal{H}\|^2$ represents the squared norm of the mean curvature $\mathcal{H}$. Moreover, the surface $\mathbb M^2$ is termed the Wintgen ideal surface when it meets the equality condition, that is, the equality is true if and only if the ellipse of the surface's curvature in $\mathbb E^4$ forms a circle. This particular inequality was further explored and expanded independently in \cite{19} and \cite{20} for surfaces with arbitrary co-dimension $n$ in the real space form $\bar M^{(n+2)}(c)$ as $$K+|K^\perp| \leq \|\mathcal{H}\|^2 +c.$$

Furthermore, B.-Y. Chen broadened the scope of the Wintgen inequality to include surfaces in pseudo-Euclidean 4-spaces $\mathbb E^4_2$ with a neutral metric, as documented in \cite{21,22}. 

In 1999, a new conjecture related to the Wintgen inequality was proposed for general Riemannian submanifolds within real space forms, which became known as the DDVV conjecture, described in \cite{23}. It was disclosed that for any submanifold $M^n$ in a real space form $\bar M^n+m (c)$, the inequality $$\rho+\rho^\perp \leq \|\mathcal{H}\|^2 +c$$ holds, where $\rho$ represents the normalized scalar curvature, and $\rho^\perp$ indicates the normalized normal scalar curvature of $M$. This inequality, also known as the generalized Wintgen inequality or the conjecture of normal scalar curvature, was independently verified by Ge and Tang \cite{24} and Lu \cite{25}.
Recent advancements in Wintgen inequalities for submanifolds within metallic Riemannian manifolds have been made. Choudhary et al. \cite{26} explored these extended Wintgen inequalities specifically for slant submanifolds in metallic Riemannian space forms equipped with semi-symmetric metric connections.

\begin{theorem} {\rm \cite{26}} \label{T:7.6}
Assume $M^n$ is a $\theta$ slant submanifold in a locally metallic space form $(\bar{M} = M_p (c_p) \times M_q (c_q),g,\varphi)$ equipped with the semi-symmetric metric connection. Then we have
\begin{equation} \label{15}
\begin{aligned}
\rho_\xi \leq & \|\mathcal{H}\|^2 - 2\rho + (c_p + c_q) \frac{p}{(p^2 + 4q)} \Bigg\{p^2 + 2q \\
& + \frac{2}{n(n-1)} [\operatorname{tr}^2 \varphi - (p \cdot \operatorname{tr} T + nq) \cos^2 {\theta}] - \frac{2p}{n} \operatorname{tr} \varphi\Bigg\} \\
& + \frac{1}{n\sqrt{p^2 + 4q}} (c_p - c_q) (2 \operatorname{tr} \varphi - np) - 2(n-1) \operatorname{tr} (\alpha)
\end{aligned}
\end{equation}
where $\alpha$ being $(0, 2)$-tensor. Moreover, the equality case of \eqref{15} holds identically if and only if, for the orthonormal frame $(E_1, \cdots, E_n, E_{n+1},\cdots, E_m ),$ the shape operator $A$ satisfy
$$
\begin{aligned}
& A_{n+1}=\left(\begin{array}{cccccc}
a & d & 0 & \ldots & 0 & 0 \\
d & a & 0 & \ldots & 0 & 0 \\
0 & 0 & a & \ldots & 0 & 0 \\
\vdots & \vdots & \vdots & \ddots & \vdots & \vdots \\
0 & 0 & 0 & \ldots & a & 0 \\
0 & 0 & 0 & \ldots & 0 & a
\end{array}\right), 
\;\;
 A_{n+2}=\left(\begin{array}{cccccc}
b+d & 0 & 0 & \ldots & 0 & 0 \\
0 & b-d & 0 & \ldots & 0 & 0 \\
0 & 0 & b & \ldots & 0 & 0 \\
\vdots & \vdots & \vdots & \ddots & \vdots & \vdots \\
0 & 0 & 0 & \ldots & b & 0 \\
0 & 0 & 0 & \ldots & 0 & b
\end{array}\right), \\
& A_{n+3}=\left(\begin{array}{cccccc}
c & 0 & 0 & \ldots & 0 & 0 \\
0 & c & 0 & \ldots & 0 & 0 \\
0 & 0 & c & \ldots & 0 & 0 \\
\vdots & \vdots & \vdots & \ddots & \vdots & \vdots \\
0 & 0 & 0 & \ldots & c & 0 \\
0 & 0 & 0 & \ldots & 0 & c
\end{array}\right), \quad A_{n+4}=\cdots=A_m=0,
\end{aligned}
$$
where $a, b, c$, and $d$ are smooth functions on $M$.
\end{theorem}

The repercussions of the Theorem \ref{T:7.6} are also discussed by Choudhary et al. in \cite{26}.

\begin{corollary} {\rm \cite{26}} \label{C:7.7}
Consider the submanifold $M^n$ immersed in a locally
metallic space form $(\bar{M} = M_p (c_p) \times M_q (c_q),g,\varphi)$ equipped with SSMC. Then, inequality \eqref{15} takes the following forms:

{\rm (i)} If $M$ is invariant 
\begin{equation} \label{16}
\begin{aligned}
\rho_\xi \leq & \|\mathcal{H}\|^2 - 2\rho + (c_p + c_q) \frac{p}{(p^2 + 4q)} \{p^2 + 2q \\
& + \frac{2}{n(n-1)} [\operatorname{tr}^2 \varphi - (p \cdot \operatorname{tr} T + nq)] - \frac{2p}{n} \operatorname{tr} \varphi\} \\
& + \frac{1}{n\sqrt{p^2 + 4q}} (c_p - c_q) (2 \operatorname{tr} \varphi - np) - 2(n-1) \operatorname{tr} (\alpha).
\end{aligned}
\end{equation}

{\rm (ii)} If $M$ is anti-invariant
\begin{equation} \label{17}
\begin{aligned}
\rho_\xi & \leq\|\mathcal{H}\|^2-2 \rho+\left(c_1+c_2\right) \frac{p}{p^2+4 q}\left[p^2+2 q+\frac{2}{n(n-1)} \operatorname{tr} \varphi^2-\frac{2 p}{n} \operatorname{tr} \varphi\right] \\
& -\frac{p}{\sqrt{p^2+4 q}}\left(c_1-c_2\right)+\frac{1}{n} \frac{1}{\sqrt{p^2+4 q}}\left(c_1-c_2\right)(2 \operatorname{tr} \varphi-n p)-2(n-1) \operatorname{tr}(\alpha) .
\end{aligned}
\end{equation}
Moreover, the equality case in \eqref{16} and \eqref{17} holds identically if and only if, for the orthonormal frame $\left\{E_1, \ldots, E_n, E_{n+1}, \ldots, E_m\right\}$, the shape operators $A$ satisfy
$$
\begin{aligned}
& A_{n+1}=\left(\begin{array}{cccccc}
a & d & 0 & \ldots & 0 & 0 \\
d & a & 0 & \ldots & 0 & 0 \\
0 & 0 & a & \ldots & 0 & 0 \\
\vdots & \vdots & \vdots & \ddots & \vdots & \vdots \\
0 & 0 & 0 & \ldots & a & 0 \\
0 & 0 & 0 & \ldots & 0 & a
\end{array}\right),\;\; A_{n+2}=\left(\begin{array}{cccccc}
b+d & 0 & 0 & \ldots & 0 & 0 \\
0 & b-d & 0 & \ldots & 0 & 0 \\
0 & 0 & b & \ldots & 0 & 0 \\
\vdots & \vdots & \vdots & \ddots & \vdots & \vdots \\
0 & 0 & 0 & \ldots & b & 0 \\
0 & 0 & 0 & \ldots & 0 & b
\end{array}\right), \\
& A_{n+3}=\left(\begin{array}{cccccc}
c & 0 & 0 & \ldots & 0 & 0 \\
0 & c & 0 & \ldots & 0 & 0 \\
0 & 0 & c & \ldots & 0 & 0 \\
\vdots & \vdots & \vdots & \ddots & \vdots & \vdots \\
0 & 0 & 0 & \ldots & c & 0 \\
0 & 0 & 0 & \ldots & 0 & c
\end{array}\right), \quad A_{n+4}=\cdots=A_m=0, \end{aligned}
$$
where $a, b, c,d$ are smooth functions on $M$.
\end{corollary}

\begin{corollary}  {\rm \cite{26}} \label{C:7.8}
Let $M^n$ be any submanifold of locally golden space form $(\bar{M} = M_p (c_p) \times M_q (c_q),g,\varphi)$ endowed with SSMC. Then we have:

{\rm (i)} If $M$ is $\theta$-slant, then
$$\begin{aligned}
\rho_\xi \leq\; &\|\mathcal{H}\|^2-2 \rho+\frac{1}{5}\left(c_1+c_2\right)\left\{3+\frac{2}{n(n-1)}\left[\operatorname{tr}^2 \varphi-(\operatorname{tr} T+n) \cos ^2 \theta\right]-\frac{2}{n} \operatorname{tr} \varphi\right\}\\&+\frac{1}{\sqrt{5} n}\left(c_1-c_2\right)(2 \operatorname{tr} \varphi-n)-2(n-1) \operatorname{tr}(\alpha) .
\end{aligned}$$

{\rm (ii)} If $M$ is invariant, then
$$\begin{aligned}
\rho_\xi  \leq\; &\|\mathcal{H}\|^2-2 \rho+\frac{1}{5}\left(c_1+c_2\right)\left\{3+\frac{2}{n(n-1)}\left(\operatorname{tr}^2 \varphi-\operatorname{tr} T-n\right)-\frac{2}{n} \operatorname{tr} \varphi\right\} \\
& +\frac{1}{\sqrt{5} n}\left(c_1-c_2\right)(2 \operatorname{tr} \varphi-n)-2(n-1) \operatorname{tr}(\alpha) .
\end{aligned}$$

{\rm (iii)} If $M$ is anti-invariant, then
$$
\begin{aligned}
\rho_\xi  \leq\; &\|\mathcal{H}\|^2-2 \rho+\frac{1}{5}\left(c_1+c_2\right)\left[3+\frac{2}{n(n-1)} \operatorname{tr}^2 \varphi-\frac{2}{n} \operatorname{tr} \varphi\right] \\
& +\frac{1}{\sqrt{5} n}\left(c_1-c_2\right)(2 \operatorname{tr} \varphi-n)-2(n-1) \operatorname{tr}(\alpha).
\end{aligned}
$$
\end{corollary}

\begin{remark}
{\rm In \cite{8}, Choudhary and Blaga discuss the generalized Wintgen inequalities for slant submanifolds in the context of metallic Riemannian space forms as well as equality instances. Moreover, the inequalities for invariant and anti-invariant submanifolds in the same ambient space are similarly found as an application of Theorem \ref{T:7.6}.}
\end{remark}

\subsection{Chen--Ricci inequality in metallic Riemannian manifolds}\label{S7.4}

The Chen-Ricci inequality is a widely recognized inequality in differential geometry that correlates a submanifold's scalar curvature to both its mean curvature and the norm of its second fundamental form.

A formula relating two geometric characteristics of a submanifold, $M$, which is embedded in a space termed $\bar M(c)$ with a constant curvature $c$, was derived in 1996 by mathematician Chen. The two characteristics are the squared mean curvature, represented by $\|\mathcal{H}\|^2$, and the Ricci curvature, represented by ``$\operatorname{Ric}$''.
Then Chen's inequality states that, for any unit vector $X$ tangent to the submanifold $M$ of $\bar M(c)$, we have $$\operatorname{Ric}(X) \leq (n-1)c+\frac{n^2}{2}\|\mathcal{H}\|^2, \;\; n=\dim M\geq 2.$$ 
 B.-Y. Chen also derived  in \cite{28} a similar inequality for Lagrangian submanifolds of a complex space form. Since its discovery, this inequality has garnered great interest from geometers worldwide. The Chen--Ricci inequality for isotropic submanifolds in locally metallic product space forms was proven by Li et al. in \cite{29}. They also ascertained the circumstances in which the inequality turns into equality.

\begin{theorem} \label{T:7.9} {\rm \cite{29}}
Consider an $n$-dimensional isotropic submanifold $M$ of an $m$-dimen\-sional locally metallic product space form $(\bar{M} = M_p (c_p) \times M_q (c_q),g,\varphi)$. Then we have:

{\rm (i)} For every unit vector $X \in T_p M$
\begin{equation} \label{18}
\begin{aligned}
    \operatorname{Ric}(X)  \leq\; & \frac{n^2}{4} \|\mathcal{H}\|^2 +\frac{1}{4} (c_p +c_q)(n-1)\left(1+\frac{p^2}{p^2 +4q}\right) \\
    & \pm \frac{1}{2}(c_p -c_q )(n-1)\frac{p}{\sqrt{p^2 +4q}}.
\end{aligned}
\end{equation}

{\rm (ii)} If $\mathcal{H}(p)=0,$ the equality case of \eqref{18} is satisfied by a unit tangent vector $X$ at $p$ if and only if $X$ in $N_p$, where $N_p$ is the relative null space of a Riemannian manifold at a point $p$ in $M.$

{\rm (iii)} If $p$ is either a totally geodesic point or if $n = 2$ and $p$ is a totally umbilical point, then \eqref{18}
equality case is true for all unit tangent vectors at $p$.
\end{theorem}

\begin{corollary} \label{C:7.10} {\rm \cite{29}}
Let $(\bar{M} = M_p (c_p) \times M_q (c_q),g,\varphi)$ be an $m$-dimensional locally metallic product space form and $M$ is n-dimensional isotropic submanifold. Then

{\rm (i)} For each unit vector $X \in T_p M$, we have
\begin{equation} \label{19}
\operatorname{Ric}(X) \leq \frac{n^2}{4} \|\mathcal{H}\|^2 +(n-1) \left[\frac{3}{10}(c_p +c_q ) \pm \frac{1}{\sqrt{5}}(c_p -c_q )\right].
\end{equation}

{\rm (ii)} If $\mathcal{H}(p) = 0,$ the equality case of \eqref{19} is satisfied by a unit tangent vector $X$ at $p$ if and only if $X \in N_p$.

{\rm (iii)} If $p$ is either a totally geodesic point or if $n = 2$ and $p$ is a totally umbilical point, then \eqref{19}
equality case is true for all unit tangent vectors at $p$.
\end{corollary}

\begin{remark}
{\rm Some inequality cases for the second fundamental form of a bi-warped product submanifold in locally nearly metallic Riemannian manifold are discussed in \cite{w5}.}
\end{remark}

\section{Inequality in Metallic-Like Statistical Manifolds}\label{S8}

\subsection{Chen-type inequality in metallic-like statistical manifolds}\label{S8.1}

In \cite{B21}, Bahadir examined the first Chen inequality as well as a Chen's inequality for the $\delta (2, 2)$ invariant for statistical submanifolds in metallic-like statistical manifolds and obtained the following results.

\begin{theorem} \label{T:8.1} {\rm \cite{B21}}
Let  $(\bar{M} = M_p (c_p) \times M_q (c_q),g,\varphi)$ be  a metallic-like statistical manifold of dimension $m$ and  $M$ be an $n$-dimensional  statistical submanifold of  $\bar{M}$. Then we have
\begin{equation*}
\begin{aligned}
&(\tau-K(\pi)) - \left(\tau_0 - K_0(\pi)\right) \geq  2\hat{K}_0(\pi) - 2\hat{\tau}_0
\\&\hskip.3in + \frac{1}{4} \frac{(c_p+c_q)(n^2-n)}{p^2+4q} \left\{1+p^2+\frac{4}{n^2-n} \left[\operatorname{tr}^2(\varphi) - \|\mathcal{T}\|^2 - \frac{4p}{n} \operatorname{tr}(\varphi)\right]\right. \\
&\hskip.3in  \pm 2 \sqrt{p^2+4q} \left(2\operatorname{tr}(\varphi) - np\right)\bigg\} - \frac{1}{4}(c_p+c_q)- \frac{c_p+c_q}{4(p^2+4q)} \left[4\Psi(\pi) - 4\Omega^2(\pi) + p^2\right] \\
&\hskip.3in  - 2p \operatorname{tr}\left(\varphi_{\left.\right|_\pi}\right) \pm \frac{1}{2} \frac{c_q-c_p}{\sqrt{p^2+4q}}\left[\operatorname{tr}\left(\varphi_{\left.\right|_\pi}\right)-p\right] - \frac{n^2(n-2)}{4(n-1)} \left[\|\mathcal{H}\|^2 + \left\|\mathcal{H}^*\right\|^2\right] ,
\end{aligned}
\end{equation*}
where $\pi \in T_p M$ is a plane section and $\Psi (\pi)= g(\varphi X, X) g(\varphi Y, Y), \Omega(\pi)=g^2(\varphi X, Y)$ for any orthonormal vectors $X, Y$ spanning $\pi$ and $\mathcal{H^*}$ refers to a mean curvature vector adapted to the statistical structure of the manifold. $\mathcal{T}X$ is the tangential component of $\varphi X$ and $\hat{K}_0, \hat{\tau}_0$ are the sectional curvature and scalar curvature according to the main statistical manifold, respectively.

Moreover, the equalities holds for any $\gamma \in\{n+1, \ldots, m\}$ if and only if
$$
\begin{gathered}
h_{11}^\gamma+h_{22}^\gamma=h_{33}^\gamma=\cdots=h_{n n}^\gamma, \\
h_{11}^{* \gamma}+h_{22}^{* \gamma}=h_{33}^{* \gamma}=\cdots=h_{n n}^{* \gamma}, \\
h_{i j}^\gamma=h_{i j}^{* \gamma}=0, i \neq j,(i, j) \neq(1,2),(2,1), 1 \leq i<j \leq n,
\end{gathered}
$$
where $h^*$ represents the second fundamental form associated with the dual connection or the statistical analogue in the context of a statistical manifold.
\end{theorem}

\begin{corollary}  \label{C:8.2} {\rm \cite{B21}}
Let $M^n$ be a totally real statistical submanifold of a metallic-like statistical manifold $(\bar{M} = M_p (c_p) \times M_q (c_q),g,\varphi)$. Then we have  
$$\begin{aligned}
&(\tau-K(\pi))-\left(\tau_0-K_0(\pi)\right) \geq 2 \hat{K}_0(\pi)-2 \hat{\tau}_0 
\\&\hskip.2in  + \frac{1}{4} \frac{\left(c_p+c_q\right)\left(n^2-n\right)}{p^2+4 q}\left\{1+p^2 \mp 2 n p \sqrt{p^2+4 q}\right\} \\
&\hskip.2in  -\frac{1}{4}\left(c_p+c_q\right)-\frac{c_p+c_q}{4\left(p^2+4 q\right)} p^2 \pm \frac{1}{2} \frac{c_p-c_q}{\sqrt{p^2+4 q}} p-\frac{n^2(n-2)}{4(n-1)}\left[\|\mathcal{H}\|^2+\left\|\mathcal{H}^*\right\|^2\right].
\end{aligned}$$
Moreover, the equalities holds for any $\gamma \in\{n+1, \ldots, m\}$ if and only if.
$$
\begin{gathered}
h_{11}^\gamma+h_{22}^\gamma=h_{33}^\gamma=\cdots=h_{n n}^\gamma, \\
h_{11}^{* \gamma}+h_{22}^{* \gamma}=h_{33}^{* \gamma}=\cdots=h_{n n}^{* \gamma}, \\
h_{i j}^\gamma=h_{i j}^{* \gamma}=0, i \neq j,(i, j) \neq(1,2),(2,1), 1 \leq i<j \leq n .
\end{gathered}
$$
\end{corollary}

\begin{corollary}  \label{C:8.3} {\rm \cite{B21}}
Let $M^n$ be a totally real statistical submanifold of a metallic-like statistical manifold $(\bar{M} = M_p (c_p) \times M_q (c_q),g,\varphi)$. If there exists a point $p \in M$ and  a plane $\pi \subset T_p M$ such that  
$$
\begin{aligned}
\tau-K(\pi)< &\; \tau_0-K_0(\pi)+\frac{1}{4} \frac{\left(c_p+c_q\right)\left(n^2-n\right)}{p^2+4 q}\left\{1+p^2 \mp 2 n p \sqrt{p^2+4 q}\right\} \\
& -\frac{1}{4}\left(c_p+c_q\right)-\frac{c_p+c_q}{4\left(p^2+4 q\right)} p^2 \pm \frac{1}{2} \frac{c_p-c_q}{\sqrt{p^2+4 q}} p+2\left(\hat{K}_0(\pi)-\hat{\tau}_0\right) .
\end{aligned}
$$
Then $M$ is non-minimal, i.e., $\mathcal{H} \neq 0$ or $\mathcal{H}^* \neq 0$.
\end{corollary}

\subsubsection{$\delta(2,2)$ Chen  inequality  in metallic-like statistical manifolds}\label{S8.1.1}

Let $p \in M$, and let $\pi_1$ and $\pi_2$ be mutually orthogonal planes within $T_p M$, spanned by $\{E_1, E_2\}$ and $\{E_3, E_4\}$, respectively. Additionally, consider $\{E_1, E_2, \ldots, E_n\}$ and $\{E_{n+1}, \ldots, E_m\}$ as the orthonormal basis of $T_p M$ and $T_p M^{\perp}$, respectively. 

The inequality, representing the $\delta(2,2)$ Chen inequality for a statistical submanifold in a metallic-like statistical manifold, is given by:

\begin{theorem} \label{T:8.4} {\rm \cite{B21}}
Let $M^n$ be a statistical submanifold of a metallic-like statistical manifold $(\bar{M} = M_p (c_p) \times M_q (c_q),g,\varphi)$. Then we have  
$$
\begin{aligned}
(\tau & -K(\pi_1)  -K(\pi_2))-(\tau_0-K_0(\pi_1)-K_0(\pi_2)) 
\\&\geq  \frac{\left(c_p+c_q\right)(n^2- n)}{4(p^2+4 q)}\Bigg[1+p^2+\frac{4 }{n^2 -n} \Big\{{\rm tr}^2(\varphi)
-\|\mathcal{T}\|^2-\frac{4p}{n} {\rm tr}(\varphi)\Big\}
\\&\hskip .2in \pm 2 \sqrt{p^2+4 q}(2 {\rm tr}(\varphi)-n p)\Bigg]
 -\frac{1}{2}(c_p+c_q)-\frac{c_p+c_q}{4(p^2+4 q)}\Big[4(\Psi (\pi_1)+\Psi (\pi_2)) \\
&\hskip .2in -4(\Omega^2(\pi_1)+\Omega^2(\pi_2))+p^2-2 p ({\rm tr}(\varphi_{|_{\pi_1}})+{\rm tr}(\varphi_{\mid \pi_2}))\Big] 
\\&\hskip .2in\pm \frac{1}{2} \frac{c_q-c_p}{\sqrt{p^2+4 q}} \big[{\rm tr}(\varphi_{|_{\pi_1}})+ {\rm tr}(\varphi_{|_{\pi_2}})-2 p\big] 
\\& \hskip .2in -\frac{n^2(n-2)}{4(n-1)}\left[\|\mathcal{H}\|^2+\left\|\mathcal{H}^*\right\|^2\right]-2\left[\hat{\tau}_0-\hat{K}_0(\pi_1-\hat{K}_0(\pi_2)\right].
\end{aligned}
$$
Moreover, the equalities holds for any $\gamma \in\{n+1, \ldots, m\}$ if and only if
$$
\begin{gathered}
h_{11}^\gamma+h_{22}^\gamma=h_{33}^\gamma+h_{44}^\gamma=h_{55}^\gamma \cdots=h_{n n}^\gamma, \\
h_{11}^{* \gamma}+h_{22}^{* \gamma}=h_{33}^{* \gamma}+h_{44}^{* \gamma}=h_{55}^{* \gamma} \cdots=h_{n n}^{* \gamma}, \\
h_{i j}^\gamma=h_{i j}^{* \gamma}=0, i \neq j,(i, j) \neq(1,2),(2,1),(3,4),(4,3), 1 \leq i<j \leq n .
\end{gathered}
$$
\end{theorem}

\begin{corollary} \label{C:8.5} {\rm \cite{B21}}
Let $M^n$ be a statistical submanifold in a metallic-like statistical manifold $\bar M$ of dimension $m$.  Then we have
 $$
\begin{aligned}
(\tau-& K(\pi_1)  -K(\pi_2))-(\tau_0-K_0(\pi_1)-K_0(\pi_2)) \\&\geq \frac{1}{4} \frac{(c_p+c_q)(n^2-n)}{p^2+4 q}\big\{1+p^2 \mp 2 n p \sqrt{p^2+4 q}\big\} \\
& -\frac{1}{2}(c_p+c_q)-\frac{c_p+c_q}{4(p^2+4 q)}\left[p^2\right] \pm \frac{1}{2} \frac{c_p-c_q}{\sqrt{p^2+4 q}}[2 p]-\frac{n^2(n-2)}{4(n-1)}\big[\|\mathcal{H}\|^2+\left[\|\mathcal{H}^*\right\|^2\big]
\\& -2\big[\hat{\tau}_0-\hat{K}_0\left(\pi_1\right)-\hat{K}_0\left(\pi_2\right)\big] .
\end{aligned}
$$
Moreover, the equalities holds for any $\gamma \in\{n+1, \ldots, m\}$ if and only if
$$
\begin{gathered}
h_{11}^\gamma+h_{22}^\gamma=h_{33}^\gamma+h_{44}^\gamma=h_{55}^\gamma \cdots=h_{n n}^\gamma, \\
h_{11}^{* \gamma}+h_{22}^{* \gamma}=h_{33}^{* \gamma}+h_{44}^{* \gamma}=h_{55}^{* \gamma} \cdots=h_{n n}^{* \gamma}, \\
h_{i j}^\gamma=h_{i j}^{* \gamma}=0, i \neq j,(i, j) \neq(1,2),(2,1),(3,4),(4,3), 1 \leq i<j \leq n .
\end{gathered}
$$    
\end{corollary}

\begin{corollary} \label{C:8.6} {\rm \cite{B21}}
A metallic-like statistical manifold $\bar M$ of dimension m has a totally real statistical submanifold denoted by $M^n$. Given a point $p \in M$ and mutually orthogonal planes $\pi_1, \pi_2 \subset T_p M$, such that
$$
\begin{aligned}
\tau-K\left(\pi_1\right)\,&-K\left(\pi_2\right)<  \tau_0-K_0\left(\pi_1\right)-K_0\left(\pi_2\right)
\\&+\frac{1}{4} \frac{\left(c_p+c_q\right)\left(n^2-n\right)}{p^2+4 q}\left\{1+p^2 \mp 2 n p \sqrt{p^2+4 q}\right\} \\
& -\frac{1}{2}\left(c_p+c_q\right)-\frac{c_p+c_q}{4\left(p^2+4 q\right)}\left[p^2\right] \pm \frac{1}{2} \frac{c_p-c_q}{\sqrt{p^2+4 q}}[2 p]
\\&-2\left[\hat{\tau}_0-\hat{K}_0\left(\pi_1\right)-\hat{K}_0\left(\pi_2\right)\right] .
\end{aligned}
$$
Then $M$ is non-minimal, i.e., $\mathcal{H} \neq 0$ or $\mathcal{H}^* \neq 0$. 
\end{corollary}

\section{Some More Structures on Metallic Riemannian \\Manifolds}\label{S9}

\subsection{Integrability of metallic Riemannian structures}\label{S9.1}

Considering the relationship between metallic structures and almost-product structures, the $\Phi$-operator method from the theory of almost-product structures can be applied to metallic structures. An integrability condition and curvature properties for these structures utilizing a $\ Phi$ operator applied to pure tensor fields are presented in \cite{in} by Gezer and Karaman.

It is well-established that a polynomial structure $F$ is integrable if and only if a torsion-free linear connection $\nabla$ can be introduced such that the structure tensor $F$ remains covariantly constant \cite{vanzura}. Utilizing the Tachibana operator, the authors in \cite{in} provided an alternative condition for the integrability of a metallic Riemannian structure.

\begin{theorem} \label{T:9.1} {\rm \cite{in}}
Consider a metallic Riemannian manifold $\bar{M}$ with a metallic structure $\varphi$ and a Riemannian metric $g$. Then: 

{\rm (i)} The structure $\varphi$ is integrable if $\Phi_{\varphi} g = 0$. 

{\rm (ii)} The condition $\Phi_{\varphi} g = 0$ is equivalent to $\bar{\nabla} \varphi = 0$, where $\bar{\nabla}$ denotes the Levi-Civita connection of $g$.
\end{theorem}

\begin{proposition}  \label{P:9.1} {\rm \cite{in}}
Consider a metallic Riemannian manifold $\bar M$ endowed with a metallic structure $\varphi$ and a Riemannian metric $g$. The manifold $\bar M$ is locally decomposable as a metallic Riemannian manifold if and only if $\Phi_{F_\pm}g=0$, where $F_\pm$ denote the almost product structures related to $\varphi$.
\end{proposition}

\subsubsection{Twin metallic Riemannian metric}\label{9.1.1}

Consider a metallic Riemannian manifold $\bar M$ with a metallic structure $\varphi$ and a Riemannian metric $g$. The twin metallic metric is defined as \cite{in}: $$G(X, Y)=g(\varphi X, Y)$$ for any vector fields $X$ and $Y$ on $\bar M$. It can be easily shown that $G$ is pure concerning $\varphi$. When the $\Phi_\varphi$-operator is applied to the metric $G$, standard computations yield $$(\Phi_\varphi G)(X, Y, Z)= (\Phi_\varphi g)(X, \varphi Y, Z) + g(\mathcal{N_\varphi} (X, Y), Z).$$

\begin{proposition}  \label{P:9.2} {\rm \cite{in}}
Let $\bar M$ be a metallic Riemannian manifold equipped with a metallic structure $\varphi$ and a
Riemannian metric $g$. Then $\Phi_\varphi g= 0$ is equivalent to $\Phi_\varphi G= 0$ if $\mathcal{N_\varphi}= 0,$ where $\mathcal{N_\varphi}$ is Nijenhuis tensor constructed from $\varphi$.
\end{proposition}

\begin{theorem}  \label{T:9.2} {\rm \cite{in}}
Let $\bar M$ be a metallic Riemannian manifold equipped with a metallic structure $\varphi$ and a Riemannian metric $g$. The Riemannian curvature tensor field is a $\varphi$-tensor field.
\end{theorem}

Results on Riemannian curvature tensor fields of the locally decomposable metallic Riemannian manifold are as follows:

\begin{proposition}  \label{P:9.3} {\rm \cite{in}}
Let $\bar M$ be a metallic Riemannian manifold equipped with a metallic structure $\varphi$ and a Riemannian metric $g$. The Riemannian curvature tensor field is a $\varphi$-tensor field.
\end{proposition}

\begin{proposition}  \label{P:9.4} {\rm \cite{in}}
Let $\bar M$ be a metallic Riemannian manifold equipped with a metallic structure $\varphi$ and a Riemannian metric $g$. The Riemannian curvature tensor field is a decomposable tensor field. 
\end{proposition}

\begin{remark}
{\rm Results on metallic structures with conformal metrics are also obtained in \cite{in}.}
\end{remark}

Manea conducted a study in \cite{some} on the integrability condition of metallic structures utilizing a Codazzi-type equation and a mixed twin metric associated with them.

Consider an $n$-dimensional Riemannian manifold $(\bar M,g).$ A $(p, q)$-metallic structure on $\bar M$ is a polynomial structure of second degree
given by a $(1, 1)$-tensor field $\varphi$ which satisfies the equation,
$$\varphi^2 -p \varphi-qI=0,$$ where $I$ is the identity on the vector fields space $\Gamma(T\bar M)$ and $p, q$ are fixed integers such that the equation $$x^2 -px-q = 0$$ has a positive irrational root
$\sigma_{p,q}$ \cite{some}.

Motivated by the Codazzi-type equation, the subsequent requirement for the metallic structure's integrability is as follows:

\begin{proposition}  \label{P:9.5} {\rm \cite{some}}
Let $(\bar M, g, \varphi)$ be a $(p, q)$-metallic manifold. The metallic structure $\varphi$ is integrable if it satisfies a Codazzi type equation, $$(d^{\bar \nabla} \varphi )(X,Y):=(\bar \nabla_X \varphi)Y-(\bar \nabla_Y \varphi)X=0,$$ for all $X,Y \in \Gamma(T\bar M),$ where $d^{\bar \nabla}$ is the exterior covariant derivative.
\end{proposition}

\begin{proposition} \label{P:9.6} {\rm \cite{some}}
For a Riemannian metallic manifold $(\bar M, g, \varphi)$ with the almost product $F$, any metallic structure induced by $F$ on $\bar M$ is integrable if
and only if the initial metallic structure is integrable.
\end{proposition}

\begin{remark}
{\rm The authors of \cite{some} have studied the metallic structures in foliated manifolds and proved that a Riemannian foliated manifold has a natural $(p, q)$-metallic structure and every leaf of the foliation carries a metallic structure that is always integrable.}    
\end{remark}

\subsection{Metallic semi-Riemannian manifolds}\label{S9.2}

The concept of a metallic Riemannian manifold can be extended to a metallic pseudo-Riemannian manifold.

\begin{definition}\cite{2018}
Let $(\bar M, g)$ be a pseudo-Riemannian manifold and let $\varphi$ be a $g$-symmetric $(1, 1)$-tensor field on $\bar M$ such that $\varphi^2 = p\varphi + qI$.
Then the pair $(\varphi, g)$ is referred to as a metallic pseudo-Riemannian structure on $\bar M$ and $(\bar M, \varphi, g)$ is known as a {\it metallic pseudo-Riemannian manifold}.
\end{definition}

\begin{definition}\cite{2018}
A linear connection $\nabla$ on $\bar M$ is termed a $\varphi$-connection if $\varphi$ is covariantly constant with respect to $\nabla$, that is, $\nabla \varphi = 0.$
\end{definition}

\begin{lemma} \label{L:9.1} {\rm \cite{2018}}
If $(\bar M, \varphi, g)$ is a locally metallic pseudo-Riemannian manifold, then $\varphi$ is integrable.
\end{lemma}

\begin{definition}\cite{2018}
A {\it trivial metallic structure} is defined as $\varphi:=\sigma I,$ where $\sigma=\frac{p\pm\sqrt{p^2 +4q}}{2}$ and $p^2 +4q\geq 0.$
\end{definition}

\begin{definition}\cite{2018}
A {\it nearly locally metallic pseudo-Riemannian manifold} $(\bar M, \varphi, g)$ is one where the Levi-Civita connection $\bar \nabla$ with respect to $g$ satisfies the condition $$(\bar \nabla_X \varphi)Y + (\bar \nabla_Y \varphi)X = 0$$ for any $X, Y \in \Gamma(TM).$
\end{definition}

\begin{proposition} \label{P:9.7} {\rm \cite{2018}}
A nearly locally metallic pseudo-Riemannian manifold $(\bar M, \varphi, g)$, where $\varphi^2 =p\varphi +qI$ and $p^2 +4q>0,$ is referred to as a locally metallic pseudo-Riemannian manifold if and only if $\varphi$ is integrable.
\end{proposition}

\subsubsection{Metallic natural connection}\label{S9.2.1}

Consider $(\bar M, \varphi, g)$ to be a metallic pseudo-Riemannian manifold where $\varphi^2 = p\varphi + qI$ and $p^2 + 4q \neq 0$. Let $\bar \nabla$ represent the Levi-Civita connection of $g$, and let $\nabla$ be linear connection defined by \cite{2018}:
\begin{equation}\label{mnc}
\nabla := \bar \nabla + \frac{2}{p^2 + 4q}\varphi (\bar \nabla \varphi) - \frac{p}{p^2 + 4q}(\bar \nabla \varphi).
\end{equation}
Consequently, we have $\nabla \varphi = 0$ and $\nabla g = 0.$ The linear connection $\nabla$ defined by \eqref{mnc} is referred to as the {\it metallic natural connection} of $(\bar M, g, \varphi)$. Torsion $T^\nabla$ of the metallic natural connection is expressed as:
$$T^\nabla (X, Y) = \frac{1}{p^2 + 4q}\{(2\varphi - pI)(\bar \nabla_X \varphi Y - \bar \nabla_Y \varphi X) - (p\varphi + 2qI)[X, Y]\}$$ for any $X, Y \in \Gamma(T\bar M).$

\begin{proposition}  \label{P:9.8} {\rm \cite{2018}}
Consider a metallic pseudo-Riemannian manifold $(\bar M, \varphi, g)$ where $\varphi^2 = p\varphi + qI$ and $p^2 + 4q \neq 0$. Torsion $T^\nabla$ of the natural connection $\nabla$ then satisfies the equation: $$T^\nabla (\varphi X, Y) + T^\nabla (X, \varphi Y) - pT^\nabla (X, Y) = (2\varphi - pI)\mathcal{N_\varphi}(X, Y),$$  $\forall X, Y \in \Gamma(T\bar M)$. Specifically, if $\varphi$ is integrable, then $$T^\nabla (\varphi X, Y) + T^\nabla (X, \varphi Y) = pT^\nabla (X, Y).$$
\end{proposition}

\subsubsection{Metallic Norden structures}\label{S9.2.2}

Note that a {\it Norden manifold} $(\bar M, \varphi, g)$ is an almost complex manifold $(\bar M, \varphi)$ with a neutral pseudo-Riemannian metric $g$ such that $g(\varphi X, Y ) = g(X, \varphi Y )$  for  $X, Y \in \Gamma(T\bar M)$.  

\begin{proposition}  \label{P:9.9} {\rm \cite{2018}}
Given a Norden manifold $(\bar M, \varphi, g)$, for all real numbers $a, b,$ $\varphi_{a,b} := a\varphi + bI$  are metallic pseudo-Riemannian structures on $\bar M$.
\end{proposition}

\begin{proposition}  \label{P:9.10} {\rm \cite{2018}}
Assume that $a \neq 0.$ Then:

{\rm (i)} $\varphi_{a,b}$ is integrable if and only if $\varphi$ is integrable,

{\rm (ii)} $\varphi_{a,b}$ is locally metallic if and only if $\varphi$ is K\"ahler,

{\rm (iii)} $\varphi_{a,b}$ is nearly locally metallic if and only if $\varphi$ is nearly K\"ahler.
\end{proposition}
On the other hand, we get the subsequent outcome:

\begin{proposition}  \label{P:9.11} {\rm \cite{2018}}
Consider the metallic pseudo-Riemannian manifold $(\bar M, \varphi, g)$ where $\varphi^2 = p\varphi + qI$ and $p^2 + 4q < 0$. Then 
$$\varphi_\pm =\pm {\small \Bigg(\frac{2}{\sqrt{-p^2-4q}}\varphi-\frac{p}{\sqrt{-p^2-4q}}I\Bigg)}$$ are Norden structures on $\bar M$ and $\varphi = a\varphi_\pm  + bI$ with $a=\pm \big({2}/{\sqrt{-p^2-4q}}\big)^{-1}$ and $b=-\frac{p}{2}$.
\end{proposition}

\begin{definition} \cite{2018}
Let $(\bar M, \varphi, g)$ be a metallic pseudo-Riemannian manifold such that $\varphi^2 = p\varphi + qI$ and $p^2 + 4q < 0$. Then $\varphi$ is known as the metallic Norden structure on $\bar M$ and
$(\bar M, \varphi, g)$ is defined as the {\it metallic Norden manifold}.
\end{definition}

\begin{remark}
{\rm In \cite{2018}, Blaga and Nannicini also discussed generalized metallic pseudo-Riemannian structures, generalized metallic natural connection, and metallic pseudo-Riemannian structures on tangent and cotangent bundles. Furthermore, Blaga and Nannicini have studied the harmonic metallic structure in \cite{harmonic}.}    
\end{remark}

\begin{remark}
{\rm Some results on metallic-like structures and metallic-like maps are obtained in \cite{2023}.}
\end{remark}

\subsection{Curvature tensors of metallic semi-Riemannian manifolds}\label{S9.3}

In \cite{2019}, Blaga and Nannicini examine the properties of curvature tensors and the concept of $\varphi$-sectional and $\varphi$-bisectional curvature of a metallic pseudo-Riemannian manifold $(\bar M, \varphi, g)$.

For non-degenerate plane sections, an analogue of the holomorphic sectional curvature will be defined as follows on a metallic pseudo-Riemannian manifold $(\bar M, \varphi, g)$.

\begin{definition} \cite{2019}
Consider $X$ to be a nonzero tangent vector field and $\pi_{X,\varphi X}$ to be the plane generated by $X$ and $\varphi X$. If $g(X,X)g(\varphi X,\varphi X)-[g(X,\varphi X)]^2 \neq 0,$ then $\pi_{X,\varphi X}$ is known as non degenerate and the $\varphi$-sectional curvature is given as:  $$K^\varphi (X) := \frac{g(R(X,\varphi X)X,\varphi X)}{g(X,X)g(\varphi X,\varphi X)-[g(X,\varphi X)]^2} ,$$ where $R$ is the Riemannian curvature tensor of $g.$ If $K^\varphi$ is a constant, then $\bar M$ has constant $\varphi$-sectional curvature denoted by $$K^\varphi(X)=R(X,\varphi X,X,\varphi X),$$ for any $X \in \Gamma(T\bar M).$
\end{definition}

\begin{proposition}  \label{P:9.12} {\rm \cite{2019}}
Assuming $(\bar M, \varphi, g)$ to be a locally metallic pseudo-Riemannian manifold, it follows that $K^\varphi (X)=0$ for any $X \in \Gamma(T\bar M)$, and the $\varphi$-sectional curvature of non-degenerate plane sections is zero.
\end{proposition}

The $\varphi$-bisectional curvature for non-degenerate plane sections on a metallic pseudo-Riemannian manifold $(M, \varphi, g)$ is defined as follows:

\begin{definition} \cite{2019}
Let $X$ and $Y$ be two linearly independent tangent vector fields, and let $\pi_{X,\varphi X}$ and $\pi_{Y,\varphi Y}$ be the planes generated by $X$ and $\varphi X$, and $Y$ and $\varphi Y$, respectively. If $\pi_{X,\varphi X}$ and $\pi_{Y,\varphi Y}$ are non-degenerate, then
the $\varphi$-bisectional curvature is defined as:
\begin{align*}
&K^\varphi (X, Y):= 
\\&{\small  \frac{g(R(X, \varphi X)Y, \varphi Y) \sqrt{(g(X,X) g(\varphi X, \varphi X) - g(X, \varphi X)^{2})(g(Y,Y) g(\varphi Y, \varphi Y) - g(Y, \varphi Y)^2)}}
{(g(X,X) g(\varphi X, \varphi X) - g(X, \varphi X)^2) (g(Y,Y) g(\varphi Y, \varphi Y) - g(Y, \varphi Y)^2)}.}
\end{align*}
If $K^\varphi$ is a constant, then $\bar M$ has constant $\varphi$-bisectional curvature given as:
$$K^\varphi (X, Y) := R(X, \varphi X, Y, \varphi Y)$$
for any $X,Y \in \Gamma(T\bar M).$
\end{definition}

\begin{proposition}  \label{P:9.13} {\rm \cite{2019}}
Given  a locally metallic pseudo-Riemannian manifold $(\bar M,\varphi ,g)$, then $K^\varphi (X, Y) = 0$ for $X,Y \in \Gamma(T\bar M)$ and the $\varphi$-bisectional curvature of nondegenerate plane sections vanishes.
\end{proposition}

\begin{theorem}  \label{T:9.3} {\rm \cite{2019}}
Consider a locally metallic pseudo-Riemannian manifold $(\bar M,\varphi,g)$ with $\varphi^2 = p\varphi + qI$ that is an RM-manifold. If $q[p^2 -(q-1)^2]\neq 0$ and $p \neq 0, q \neq -1,$ then its Riemann curvature tensor vanishes.
\end{theorem}

\vskip.2in
\noindent {Chen: Department of Mathematics, Michigan State University, East Lansing, Michigan 48824--1027, U.S.A.}
\email{chenb@msu.edu}

\vskip.2in
\noindent {Choudhary: Department of Mathematics, School of Sciences, Maulana Azad National Urdu University, Hyderabad, India}\email{majid.alichoudhary@gmail.com}

\vskip.2in
\noindent {Perween: Department of Mathematics, School of Sciences, Maulana Azad National Urdu University, Hyderabad, India}\email{afshanperween99@gmail.com}

\end{document}